\documentclass[mnsc,nonblindrev]{informs3} 

\usepackage{hyperref}
\OneAndAHalfSpacedXI 

\usepackage{bm,amsmath,tabularx}
\newcolumntype{L}[1]{>{\raggedright\arraybackslash}p{#1}}
\newcolumntype{C}[1]{>{\centering\arraybackslash}p{#1}}
\newcolumntype{R}[1]{>{\raggedleft\arraybackslash}p{#1}}

\usepackage[dvipsnames]{xcolor}
\newcommand{\beq}{\begin{equation}}
\newcommand{\eq}{\end{equation}}
\newcommand{\E}{\mathbb{E}}

\renewcommand{\P}{\mathbb{P}}

\newcommand{\das}{\stackrel{d}{=}}

\renewcommand{\vec}[1]{\mathbf{#1}}

\newcommand{\eqan}[1]{\begin{align} #1 \end{align}}
\newcommand{\e}{{\rm e}}

\newcommand{\expect}[1]{{\mathbb E}[#1]}
\newcommand{\expectp}[1]{{\mathbb E}_{\mathbb P}[#1]}
\newcommand{\prob}[1]{{\mathbb P}(#1)}
\newcommand{\indep}{\protect\mathpalette{\protect\independenT}{\perp}}
\def\independenT#1#2{\mathrel{\rlap{$#1#2$}\mkern2.5mu{#1#2}}}

\newcommand{\eps}{\varepsilon}
\newcommand{\il}{\int\limits}

\newcommand{\Dist}{\mathbb{P}}

\definecolor{darkgreen}{rgb}{0.0, 0.4, 0.13}

\definecolor{dgreen}{rgb}{0,0.7,0}
\definecolor{dred}{rgb}{0.8,0,0}
\definecolor{dblue}{rgb}{0,0,0.8}

\newcommand{\UU}{U}
\newcommand{\VV}{V}

\usepackage{natbib}
 \bibpunct[, ]{(}{)}{,}{a}{}{,}%
 \def\bibfont{\small}%
\TheoremsNumberedThrough     
\ECRepeatTheorems

\EquationsNumberedThrough    

\MANUSCRIPTNO{MS-0001-1922.65}

\usepackage{pgfplots}
\usepackage{tikz}
\usetikzlibrary{shapes,arrows,positioning}
\usetikzlibrary{automata}
\usepackage{bbm}

\pgfplotsset{compat=newest}
\pgfplotsset{plot coordinates/math parser=false}
\pgfplotsset{
    every non boxed x axis/.style={
        xtick align=center,
        enlarge x limits=true,
        x axis line style={line width=0.8pt, -latex}
},
    every boxed x axis/.style={}, enlargelimits=false
}
\pgfplotsset{
    every non boxed y axis/.style={
        ytick align=center,
        enlarge y limits=true,
        y axis line style={line width=0.8pt, -latex}
},
    every boxed y axis/.style={}, enlargelimits=false
}
\usetikzlibrary{arrows,intersections}

\renewcommand{\theARTICLETOP}{}
\usepackage{fancyhdr}
\renewcommand\headrulewidth{0pt}
\begin{document}




\RUNTITLE{MAD dispersion measure makes extremal queue analysis simple}

\TITLE{
MAD dispersion measure makes extremal\\ queue analysis simple}

\ARTICLEAUTHORS{%
\AUTHOR{Wouter van Eekelen}
\AFF{Department of Econometrics and Operations Research, Tilburg University, \EMAIL{w.j.e.c.vaneekelen@tilburguniversity.edu}} 
\AUTHOR{Dick den Hertog}
\AFF{Amsterdam Business School, University of Amsterdam, \EMAIL{d.denhertog@uva.nl}} 
\AUTHOR{Johan S.H. van Leeuwaarden}
\AFF{Department of Econometrics and Operations Research, Tilburg University, and Department of Mathematics and Computer Science, Eindhoven University of Technology,  \EMAIL{j.s.h.vanleeuwaarden@tilburguniversity.edu}}
} 

\ABSTRACT{%
A notorious problem in queueing theory is to compute 
the worst possible performance of the GI/G/1 queue under mean-dispersion constraints for the interarrival and service time distributions. We address this extremal queue problem by measuring dispersion in terms of Mean Absolute Deviation (MAD) instead of variance, making available recently developed techniques from Distributionally Robust Optimization (DRO). Combined with classical random walk theory, we obtain explicit expressions for the extremal interarrival time and service time distributions, and hence the best possible upper bounds, for all moments of the waiting time. {We also apply the DRO techniques to obtain tight lower bounds that together with the upper bounds provide robust performance intervals. We show that all bounds are computationally tractable and remain sharp, also when the mean and MAD are not known precisely, but estimated based on available data instead.}

}%



\KEYWORDS{extremal queue problem, GI/G/1 queue, random walk theory, tight bounds, distributionally robust optimization} 


\maketitle

%


\section{Introduction}


Queueing theory exists for more than a century with throughout a central role for the GI/G/1 queue with i.i.d.~interarrival times $ \{{\UU}_n \} $ distributed as $\UU$ and i.i.d.~service times $ \{{\VV}_n \} $ distributed as $\VV$. The waiting times in the GI/G/1 queue can be expressed as the maxima of a random walk with step size $ X = \VV-\UU $, the subject of an enormous literature: \cite{chung2001course,feller1971,asmussen2003}.
For all moments of the maxima (i.e.,~waiting times), general expressions are available that involve convolutions of the distribution of $X$. To use these general expressions, one thus needs to specify the precise distribution of $X$, and in the case of the GI/G/1 queue the distributions of both $\UU$ and $\VV$.

Special cases of the GI/G/1 queue can be studied with dedicated techniques for Markov chains. For instance, the M/G/1 queue with Poisson arrivals and the GI/M/1 queue with exponential services have explicit solutions that are more insightful than the general random walk results: \cite{asmussen2003,cohen1982}.
Another large, somewhat opposite branch of queueing theory concerns finding approximations and bounds.  
For the steady-state waiting time $W$ in the GI/G/1 queue, the arguably most famous upper bound for $\expect{W}$ was obtained by \cite{kingman1962some}  in terms of the first two moments of both $\UU$ and $\VV$.
While Kingman's bound is sharp in situations of heavy traffic, when $ \expect{\UU} /\expect{\VV}$ approaches 1, it leaves room for improvement for all other values of $ \expect{\UU} /\expect{\VV}$. 

In search for that sharpest possible (tight) upper bound under the first two moments constraints, foundational work was done by \cite{rolski1972some,eckberg1977sharp}, and \cite{whitt1984approximations} in the context of the GI/M/1 queue. 
\cite{whitt1984approximations} considered the GI/M/1 queue with given mean and variance of $\UU$, and showed that $\expect{W}$ is maximized when the interarrivals follow a specific two-point distribution. 
It also led to the conjecture that the overall worst case behavior (in terms of $\expect{W}$) would be caused by two-point distributions, for both $\UU$ and $\VV$. That conjecture was proved invalid by counterexamples in \cite{whitt1984approximations} when fixing either $\UU$ or $\VV$, but the conjecture remained standing for 
the case when both $\UU$ and $\VV$ are unspecified, except for their first two moments.  After that it remained silent for a while, until \cite{chen2019extremal} showed recently, for distributions with finite support, that the extremal distributions of $\UU$ and $\VV$ both have supports on at most three points.  While existence is thus proved, the exact form of the extremal three-or-fewer-points distributions can only be determined numerically, as the solution of a hard non-convex nonlinear optimization problem. Extensive numerical experiments led Chen and Whitt to conjecture that the worst case 
 is formed by two-point distributions for both $\UU$ and $\VV$, in line with the conjecture postulated several decades ago. Finding the extremal queue for given mean-variance information is therefore one of the longest standing problems in the field. That problem remains open, also after publication of the present paper.

We do consider the same problem of finding the sharpest possible bounds for GI/G/1 queue performance measures, but take a radical turn by quantifying dispersion in terms of mean absolute deviation (MAD) instead of variance. That may appear a bold decision, because MAD is hardly used in queueing theory, or random walk theory for that matter. We can only speculate about the historical reasons for  variance preference, but the random walk and GI/G/1 queue are intrinsically linked with i.i.d.~sums of random variables, and variance then enters naturally (e.g.,~variance of the sum, central limit theorem). The variance and MAD, however, are equally adequate descriptors of dispersion, and are both easily calibrated on data using basic statistical estimators. 

The MAD perspective offered in this paper departs from the variance-based formulations of the past (see \cite{rolski1972some,eckberg1977sharp,whitt1984approximations,
chen2019extremal} and the references therein), and brings to bear the rich theory of robust optimization, in particular the rapidly expanding theory of distributionally robust optimization (DRO). The exact expressions for the random walk maxima form a crucial ingredient for our proof methodology. These expressions are convex functions of 
the driving random variables, 
a prerequisite for the mean-MAD approach. Indeed, recent advances in DRO, see \cite{postek2018robust}, show that knowledge on the support, mean and MAD can lead to closed-form expressions for stochastic quantities such as the minimum and maximum expectation of a convex function.

{Using the MAD instead of the variance as dispersion measure has several important advantages for, e.g., analyzing the waiting times in
 GI/G/1 queues.
First, not only simple explicit expressions for the worst-case distributions can be obtained, but also for the best-case ones. Hence, a sharp upper bound and a sharp lower bound for the expected waiting time can be obtained. Second, our approach is for i.i.d.~sums of random variables, while existing DRO approaches have to tolerate possible dependence structures between the random variables. Third, our approach is suitable for analyzing both transient behavior and the steady state. Fourth, because of its computational tractability our approach can also be extended to many optimization variants.}

The contributions of this paper can be summarized as follows:

\begin{henumerate}
\item We suggest to use MAD instead of variance, and obtain by concise mathematical proof the worst-case three-point distribution for a rich class of extremal problems. 
This proof for MAD gives insight into why the traditional moment constraints, although a popular choice, may not necessarily yield tractable counterparts.
\item  We leverage this result to obtain tight upper and lower bounds for performance measures, including transient and steady-state queue length moments. Under mean-MAD constraints, these bounds are the sharpest possible (and thus cannot be improved). The mean-MAD approach in this paper is a new quantitative method applicable to random walks, queues and related stochastic processes. This generic approach is a computationally tractable way to analyze key performance measures of such processes.


\item We present guidelines that describe how to compute the novel tight bounds efficiently. Moreover, we demonstrate our approach when the mean and MAD are not known precisely and need to be estimated from data. Also in these more realistic settings, the bounds remain sharp.
\end{henumerate}

\noindent{\bf Outline.} The remainder of the paper is organized as follows. 
Section \ref{sec2}  presents the MAD perspective. Section \ref{sec3} discusses methods to obtain upper and lower bounds for both best and worst-case performance. Section \ref{sec4} presents a full solution of the extremal queue problem with mean-MAD constraints, and draws a comparison with the traditional mean-variance setting. We conclude in Section \ref{sec6}, also
mentioning  possibilities for follow-up research. 

\noindent{\bf Notation.} Boldfaced characters represent vectors, and $x_i$ denotes the $i$-th element of vector $\vec{x}$. For a random variable $X$, we use $X\sim \mathbb{P}\in \mathcal{P}$ to say that $X$ is a random variable with probability distribution $\mathbb{P}$ from the set of probability distributions $\mathcal{P}$.  We denote $\mathbb{E}_{\mathbb{P}}[\cdot]$ as the expectation over the probability distribution $\mathbb{P}$. When we consider $\mathbb{E}_{\mathbb{P}}[f(\vec{X})]$ with $\vec{X}=(X_1,\ldots,X_n)$, 
it is tacitly assumed that $f(\cdot)$ is
a measurable function from $\mathbb{R}^n$ to $\mathbb{R}$, and  such that $\mathbb{E}_{\mathbb{P}}[f(\vec{X})]$ exists.

\section{Extremal random walk}\label{sec2}

Consider the partial sums $S_n:=X_1+\cdots+ X_n$ ($S_0:=0$) of i.i.d.~random variables $X_1,X_2,\ldots$ distributed as $X$. The random walk $(S_n,n\geq 0)$ arises in many application domains, including queueing theory, inventory management and risk theory. If  $(S_n,n\geq 0)$ indeed models congestion, shortfall or capital position, large values of $S_n$ are of particular interest, and it is natural to consider the maxima sequence $M_n:=\max\{S_0,S_1,\ldots,S_n\}$. The random walk 
and its maxima can be studied with mathematical techniques for sums of random variables, covered in many standard texts on probability theory, e.g.,~\cite{asmussen2003,chung2001course,cohen1982,feller1971}. For the distribution and moments of $M_n$ there exist general formulas in terms of finitely many convolutions. However,  applying these exact formula requires full specification of the distribution of $X$. This paper searches for the sharpest possible bounds on 
$\expect{M_n}$ and related quantities, when only information is available on the mean and dispersion of $X$. 
We now present such bounds when the partial information consists of the mean, range and MAD of $X$.

\subsection{Extremal distribution}
Notice that $M_n$ can be expressed as $h_n(X_1,\ldots,X_n)$, with
\begin{equation}
h_n(x_1,\ldots,x_n)=\max\{0,x_1,\ldots,x_1+\cdots+x_n\}, 
\end{equation}
and the expected maximum can be expressed as $\expect{M_n}=\expect{h_n({\vec{X}})}$ with ${\vec{X}}=(X_1,\dots,X_n)$. 
For now assume that $X_1,\ldots,X_n$ are independent, but that each 
$X_i$ can have a different distribution. Assuming we only have partial information consisting of means and dispersion measures of the random variables $X_1,\ldots,X_n$, the first question we ask and answer in this paper is: What {\it extremal} distributions of  $X_i$ result in the worst-case expected maxima? Extremal distributions have been studied in many contexts, and in the literature variance
is predominantly used as the dispersion measure. Here we shall use the MAD. 
To describe all considered distributions we define an ambiguity set  that consists of all distributions of componentwise independent $\vec{X}$ with known supports, means, and MADs. 
The partial information for  $(X_1,\ldots,X_{n})$ consists of (i) $X_i$ has support  $ \text{supp}(X_i)= [a_i,b_i]$ with $-\infty < a_i \leq b_i < \infty, i=1,\ldots,n$, 
(ii) $\mathbb{E}_\mathbb{P}(X_i) = \mu_i $ and  (iii) 
$\mathbb{E}_\mathbb{P} |X_i - \mu_i| = d_i $. 
This defines the \textit{ambiguity set}
\begin{equation}
\label{eq:mu_d_conditions_multivariable}
\mathcal{P}_{(\mu,d)} = \left\{ \mathbb{P}: \ \text{supp}(X_i) \subseteq [a_i,b_i], \ \mathbb{E}_{\mathbb{P}}(X_i) = \mu_i, \ \mathbb{E}_{\mathbb{P}}\left| X_i - \mu_i \right| = d_i,\ \forall i, \ X_i \indep X_j, \ \forall i \neq j \right\},
\end{equation}
where $X_i \indep X_j, \, \forall i \neq j,$ denotes stochastic independence of the components $X_1,\dots,X_n$. In what follows, $\vec{X}$ is a  vector of random variables whose distribution $\mathbb{P}$ belongs to the set $\mathcal{P}_{(\mu,d)}$. 

As the title says, with MAD as dispersion measure, the extremal problem becomes simple. Observe that the function $h_n$ is convex in the vector 
$(x_1,\ldots,x_n)$. We can thus apply the general upper bound in \cite{BenTal1972} on the expectation of a convex function of independent random variables with mean-MAD ambiguity, which gives the following result:
\begin{theorem}\label{thm1h}
The extremal distribution that solves 
\begin{equation}\label{maxphrase2}
\max\limits_{\mathbb{P} \in \mathcal{P}_{(\mu,d)}} \expectp{h_n(\vec{X})}
\end{equation}
consists for each $X_i$ of a three-point distribution with values 
$\tau_1^{(i)} = a_i$, $\tau_2^{(i)}  = \mu_i$, $\tau_3^{(i)}  = b_i $
and probabilities
\begin{equation}
\label{eq:BenTal_probabilities_multivariable}
p_{1}^{(i)}  = \frac{d_i}{2(\mu_i - a_i)}, \quad p_{2}^{(i)}  = 1- \frac{d_i}{2(\mu_i - a_i)}  - \frac{d_i}{2(b_i - \mu_i)} , \quad p_{3}^{(i)}  = \frac{d_i}{2(b_i - \mu_i)}.
\end{equation}
\end{theorem}



 {\cite{BenTal1972} prove Theorem~\ref{thm1h} (for general convex functions) by introducing 
 a piecewise linear function on the interval $[a,b]$ that intersects the convex function in $a$, $\mu$ and $b$, and then applying the classic Jensen bound to the subintervals 
 $[a,\mu]$ and $[\mu,b]$. In the next section, we give another proof of Theorem~\ref{thm1h} that also gives insight into why using as dispersion measure MAD instead of variance makes the analysis so simple.  }


\subsection{{Novel primal-dual proof of Theorem~\ref{thm1h}}}\label{novelpr}
Our proof will crucially rely on the fact that the univariate case of Theorem~\ref{thm1h} is tractable, and can be straightforwardly extended to the multivariate case. We thus start by considering some 
 univariate measurable function $f(x)$ (with the univariate function $h_1(x_1)$ as an example) that has finite values on $[a,b]$, the support of the distribution $\Dist(x)$. 
Under mean-MAD ambiguity of one random variable $X$ we thus need to solve
\begin{equation}\label{test3}
\begin{aligned}
&\max_{\Dist(x)\geq0} &  &\int_x f(x){\rm d} \Dist(x)\\
&\text{s.t.} &      & \int_x |x-\mu|{\rm d}\Dist(x)=d, \int_x x{\rm d}\Dist(x)=\mu,\ \int_x {\rm d}\Dist(x)=1,   
\end{aligned}
\end{equation}
 a semi-infinite linear program (LP) with three equality constraints. 
{ A perhaps surprising, yet classical fact, is that the semi-infinite LP \eqref{test3} can be reduced to an equivalent finite LP that yields the same optimal value. 
Indeed, the Richter-Rogosinski Theorem (e.g.,~\cite{rogosinski1958moments,shapiro2009lectures,han2015convex})
states that there exists an extremal distribution for problem  \eqref{test3} with at most three support points. While finding these points in closed form is typically not possible (for general semi-infinite problems), we next show that this is possible for the problem at hand, by resorting to the dual problem and exploiting both the specific shape of the MAD constraint $\int_x |x-\mu|{\rm d}\Dist(x)=d$ and convexity of $f$. }



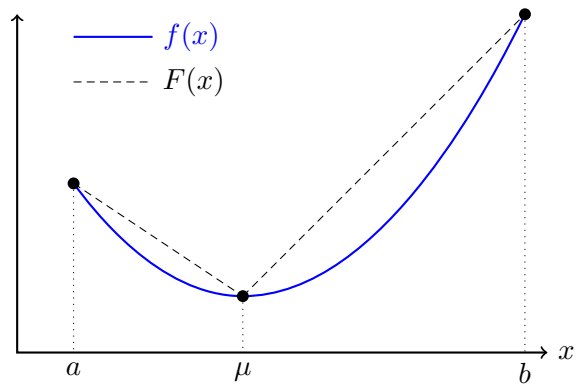
\begin{figure}[h]
\begin{center}
\begin{tikzpicture}[scale=1.5]
\draw [<->,thick] (0,3) node (yaxis) [above] {}
        |- (4.7,0) node (xaxis) [right] {$x$};
\draw[blue,thick] (0.5,1.5) parabola bend (2,.5) (4.5,3)
        node[below right] {};
        \draw[blue,thick] (0.5,2.8) -- (1.2,2.8) node[right] {$f(x)$};
        \draw[densely dashed] (0.5,2.4) -- (1.2,2.4) node[right] {$F(x)$};
        \coordinate (a) at (0.5,1.5);
        \fill[black] (a) circle (1.5pt);
        \coordinate (b) at (2,.5);
        \fill[black] (b) circle (1.5pt);
        \coordinate (c) at (4.5,3);
        \fill[black] (c) circle (1.5pt);
         \draw[densely dashed] (a) -- (b);
         \draw[densely dashed] (b) -- (c);
         \draw[dotted] (b) -- (2,0) node[below] {$\mu$};
          \draw[dotted] (a) -- (0.5,0) node[below] {$a$};
           \draw[dotted] (c) -- (4.5,0) node[below] {$b$};
\end{tikzpicture}
\caption{Some convex function $f(x)$ and its piecewise linear majorant $F(x)$.}\label{fig:major2}
\end{center}
\end{figure}

Consider the dual of \eqref{test3},
\begin{equation}\label{test4}
\begin{aligned}
&\min_{\lambda_1,\lambda_2, \lambda_3} &  &\lambda_1 d+\lambda_2 \mu+\lambda_3\\
&\text{s.t.} &      & f(x)-\lambda_1|x-\mu|-\lambda_2x-\lambda_3\leq 0, \ \forall x\in[a,b].
\end{aligned}
\end{equation}
Define $F(x)=\lambda_1|x-\mu|+\lambda_2x+\lambda_3$. Then the inequality in \eqref{test4} can be written as $f(x)\leq F(x)$, $\forall x$, i.e.~$F(x)$ majorizes $f(x)$. Note that $F(x)$ has a `kink' at $x=\mu$. Since the dual problem \eqref{test4} has three variables, the tightest majorant $F(x)$ touches $f(x)$ at three points: $x=a$, $\mu$ and $b$, as illustrated in Figure~\ref{fig:major2}.  The optimal probabilities of \eqref{test3} can now easily be obtained by solving the linear system resulting from the equations of \eqref{test3}. This is a linear system of three unknown probabilities and three equations, with a solution as stated in Theorem~\ref{thm1h}.

%
To deal with the multivariate case,
 we recursively apply the univariate result. Suppose we first apply this result to $x_1$, then the worst-case distribution is as in Theorem~\ref{thm1h}, independent of the values for $x_2,\ldots,x_n$. Moreover, the worst-case expectation  becomes a convex function in $x_2,\ldots,x_n$, since the worst-case probabilities for $x_1$ are nonnegative. Hence, we can apply the result above for the univariate case to $x_2$, etc. This completes the proof. Note that for multiperiod problems that involve multivariate optimization, such as the waiting time in the GI/G/1 queue, determining the extremal distribution for period $n$ is unaffected by all previous periods.

{



To the best of our knowledge, our proof is the first to exploit the specific shape of the kink-majorant to find an analytic solution for the semi-infinite LP. While the dual problems are often solvable as semi-definite or second-order conic programs, analytic solutions as in our case are typically hard to attain, and require special structural properties of the LP's objective function or its constraints. Notice that in the univariate case, this proof method does not require convexity of $f(x)$ and in fact could work for an arbitrary measurable function $f$. Convexity is needed, however, in the proof of Theorem~\ref{thm1h} to extend the univariate case to the multivariate case.
The proof method is of independent interest, and can for instance be applied to study the mean-MAD counterparts of the mean-variance analyses in e.g.~\cite{xin2013time,natarajan2007mean,perakis2008regret,semivariance}, and \cite{das2018heavy}.
}

\subsection{Why is MAD computationally easier than variance?}

Now that we fully grasp why and how the proof of  Theorem~\ref{thm1h} relies on the specific structural properties of the mean-MAD constraints, and in particular the univariate result seamlessly passes into the multivariate counterpart, we can also explain why the comparable challenge with mean- variance constraints becomes much more difficult if not impossible.
Observe that for the univariate case, the same proof argument works when $\sigma^2$ is given instead of $d$, i.e.,~when $|x-\mu|$ in \eqref{test3} is replaced by $(x-\mu)^2$. Hence, irrespective of whether MAD or variance is used as dispersion measure, for determining the tight upper bound of $f(x)$, it suffices to consider distributions with support on at most three points. There is however a crucial complication when extending to the multivariate case. 

To see this, observe  that when $\sigma^2$ is used as dispersion measure, the end points and kink point do {\it not} necessarily span the support of the extremal distribution. That is, upon replacing $|x-\mu|$ with $(x-\mu)^2$, the tightest majorant $F(x)$ does not necessarily touch $f(x)$ in $a$, $b$ and $\mu$. 

Hence, if the variance is used as dispersion measure, then the worst-case distribution depends on the function $f(x)$. This has severe consequences for the multivariate case, i.e., when we consider $h_n(x_1,\ldots,x_n)$. In that case, 
the worst-case distribution depends on the values of $x_2,\ldots,x_n$, and calculating (in closed form) the worst-case distribution as a function of $x_2,\dots,x_n$ seems to be impossible. Moreover, even if we would be able to derive such a worst-case distribution, substituting this distribution in the worst-case expectation would result in an extremely difficult function in $x_2,\ldots,x_n$ that is likely non-convex, and hence applying the univariate result to $x_2$ is no longer possible. 
{Our duality proof thus reveals that the complicating feature of the mean-variance framework applied to multiperiod problems is the fact that the extremal distribution in period $n$ is affected by all previous periods.}



\section{Sharpest possible bounds}\label{sec3}
A direct consequence of Theorem \ref{thm1h} is that the worst-case expectation of $h_n(\vec{X})$ is obtained by  enumerating over all $3^{n}$ permutations of outcomes $a_i,\mu_i, b_i$ of components $X_i$. 
\begin{corollary}\label{corollary1}
\begin{equation}
\label{permm4}
\max\limits_{\mathbb{P} \in \mathcal{P}_{(\mu,d)}} \expectp{h_n(\vec{X})} = \sum\limits_{\bm{\alpha} \in \{ 1,2,3 \}^{n}}  h_n(\tau_{\alpha_1}^{(1)} ,\ldots, \tau_{\alpha_{n}}^{(n)})\prod\limits_{i=1}^{n} p_{\alpha_i}^{(i)}.
\end{equation}
\end{corollary}
Thus, under the partial information contained in $\mathcal{P}_{(\mu,d)}$, \eqref{permm4} is an upper bound on $\expect{M_n}$ 
 that cannot be improved. We next specialize to the random walk setting with $X_1,X_2,\ldots$ independent and distributed as $X$, obtain representations for the tight upper bound that are computationally less cumbersome than \eqref{permm4}, and extend to all moments of the all-time maximum (when $n\to\infty$). 

\subsection{Random walk upper bounds}\label{}
We recall that
 \cite{spitzer1956} used combinatorial arguments to establish for $\expect{M_n}$ the alternative expression (which strictly requires i.i.d. increments)
\begin{equation}\label{permm3}
\expect{M_n}=\sum_{k=1}^n\frac1k \expect{S_k^+},
\end{equation}
with $x^+=\max\{0, x\}$. This can be written as $\expect{M_n}=\expect{f_n({\vec{X}})}$  with 
\begin{equation}
f_n(x_1,\ldots,x_n)=\sum_{k=1}^n \frac1k \max\{0,x_1+\ldots+x_k\}.
\end{equation}

A first usage of Spitzer's formula  \eqref{permm3}  is a considerable improvement, in terms of computational complexity, of the tight bound for $\expect{M_n}$ in \eqref{permm4}.  To state the result and for later reference, let $\Omega(\mu,d,a,b)$ denote a three-point distribution on 
the values $\{a,\mu,b\}$ with probabilities
\begin{equation}
\label{basicprobs}
p_{1}= \frac{d}{2(\mu- a)}, \quad p_{2} = 1- \frac{d}{2(\mu - a)}  - \frac{d}{2(b- \mu)} , \quad p_{3} = \frac{d}{2(b - \mu)}.
\end{equation}
Let $X_{(3)}$ denote the random variable with the extremal three-point distribution, identified in Theorem~\ref{thm1h} for the special case when $X_1,X_2,\ldots$ are i.i.d., hence $X_{(3)}\sim \Omega(\mu,d,a,b)$.

\begin{corollary}
\begin{equation}
\label{permm6}
\max\limits_{\mathbb{P} \in \mathcal{P}_{(\mu,d)}} \expectp{f_n(\vec{X})}=
\sum_{k=1}^n \frac1k \sum_{\sum_i k_i=k}\max\{0,k_1a+k_2\mu+k_3b\}\cdot \frac{k!}{k_1!k_2!k_3!}p_1^{k_1}p_2^{k_2}p_3^{k_3}.
\end{equation}
\end{corollary}
Note that for each fixed $k$, \eqref{permm6} contains a multinomial distribution with support set 
$
\{(k_1,k_2,k_3)\in\mathbb{N}^3:k_1+k_2+k_3=k\}
$
with cardinality ${k+2\choose 2}$. This implies that the sum over $k$ in \eqref{permm6} is over roughly $n^3$ terms, which is way better than the $3^n$ terms in \eqref{permm4}. 

For $\expect{X}<0$ the all-time maximum $M:=\lim_{n\to\infty}M_n$ is a proper random variable ($M_n$ converges in distribution to $M$, which will be finite with probability one if $\expect{X}<0$). 
Let $c_m(M)$ denote the $m$-th cumulant of $M$. Recall that  $c_1(M)$ is the mean,  $c_2(M)$ is the variance, and  $c_3(M)$ is the central moment $\expect{(M-\expect{M})^3}$. From general random walk theory we know that (see e.g.,~\cite{Abate1993})
\begin{equation}\label{eq: thm2}
c_m(M)=\sum_{k=1}^\infty \frac1k \expect{(S_k^+)^m}.
\end{equation}
We can now prove results similar as for $\expect{M_n}$, regarding the extremal distribution and tight upper bound. 


\begin{theorem}\label{thm1f}
Consider the random walk with generic step size $X$ contained in the ambiguity set $\mathcal{P}_{(\mu,d)}$. The tight upper bounds for all cumulants $c_m(M)$ of the all-time maximum $M$ are the cumulants of the random walk with extremal step size $X_{(3)}$. 
\end{theorem}
\proof{Proof.}
Consider the function
\begin{equation}
f_n^m(x_1,\ldots,x_n)=\sum_{k=1}^n \frac1k \left(\max\{0,x_1+\ldots+x_k\}\right)^m,
\end{equation}
which is  convex in the vector 
$(x_1,\ldots,x_n)$. Hence, for i.i.d.~increments with generic $X$, 
\begin{equation}\label{maxphrase2}
\max\limits_{\mathbb{P} \in \mathcal{P}_{(\mu,d)}} \expectp{f_n^m(\vec{X})}
\end{equation}
is solved by the extremal random variable $X_{(3)}$. This gives the bound, with $X_1^*,X_2^*,\ldots$ i.i.d.~as $X_{(3)}$, 
\begin{equation}
l_n:=\sum_{k=1}^n \frac1k \expect{(S_k^+)^m}\leq \mathbb{E} f_n^m(X_1^*,\ldots, X_n^*)=:u_n.
\end{equation}
The result follows by observing that the sequences $\{l_n\}$ and $\{u_n\}$ are both monotone, and converging to well-defined limits. 
\Halmos
\endproof

We conclude that the extremal three-point distribution for $\expect{M_n}$ in Theorem~\ref{thm1h} is also the extremal distribution for all cumulants of $M$. When calculating the associate tight upper bounds for $c_m(M)$, \eqref{eq: thm2} shows that we are confronted with an infinite summation of increasingly complex summands. Here, another line of classical random walk theory can help, which transforms such infinite sums into complex contour integrals. 

Consider the random walk with generic step size $X$. It is known that formal solutions of the distribution of $M_n$ and $M$ can 
be expressed in terms of complex contour integrals (see \cite{Abate1993,Janssen2015} for the algorithmic aspects of these contour integrals). 
Assume  that $\phi_X(s)=\expect{{\rm e}^{sX}}$ is analytic for complex $s$ in the strip $|{\rm Re}(s)|<\delta$ for some $\delta>0$. A  sufficient condition is that the moment generating function $\phi_X(s)$ is finite in a neighborhood of the origin, and hence all moments of $X$ exist. Then 
\begin{equation}\label{polgen}
\expect{{\rm e}^{-sM}}=
\exp\left\{\frac{-1}{2\pi i}\int_{\mathcal{C}}\frac{s}{u(s-u)}\log(1-\phi_X(-u)){\rm d}u\right\},
\end{equation}
where $s$ is a complex number with ${\rm Re}(s)\geq 0$, $\mathcal{C}$ is a contour to the left of, and parallel to, the imaginary axis, and to the right of any singularities of $\log(1-\phi_X(-u))$ in the left half plane. From  \eqref{polgen}
contour integral expressions
for the cumulants follow by differentiation:
\begin{equation}\label{polgen2}
c_m(M)= \frac{(-1)^m}{2\pi i}\int_{\mathcal{C}}\frac{\log(1-\phi_{X}(-u))}{u^{m+1}}{\rm d}u. 
\end{equation}



 Consider $X=X_{(3)}$ with a three-point distribution on values $\{a,b,c\}$ with probabilities $p_a,p_b,p_c$ and moment generating function
\begin{equation}\label{mgfpol}
\phi_{X_{(3)}}(s)=p_a  {\rm e}^{sa}+ p_b {\rm e}^{sb}+p_c{\rm e}^{sc}.
\end{equation}
Notice that all moments of $X_{(3)}$ exist, and hence $\phi_{X_{(3)}}(s)$  satisfies the assumption required for representation \eqref{polgen} to hold. Since $X_{(3)}$ follows the extremal three-point distribution associated with the tight upper bounds for $c_m(M)$, we obtain the following result:

\begin{corollary}Let $\phi_{X_{(3)}}(s):=\expect{{\rm e}^{sX_{(3)}}}=p_1  {\rm e}^{sa}+ p_2 {\rm e}^{s\mu}+p_3{\rm e}^{sb}.$ The tight upper bounds on $c_m(M)$ identified in {\rm Theorem~\ref{thm1f}} are given by
\begin{equation}\label{pollaczekup}
 \frac{(-1)^m}{2\pi i}\int_{\mathcal{C}}\frac{\log(1-\phi_{X_{(3)}}(-u))}{u^{m+1}}{\rm d}u, \quad m=1,2,\dots,
\end{equation}
where $\mathcal{C}$ is a contour to the left of, and parallel to, the imaginary axis, and to the right of any singularities of $\log(1-\phi_{X_{(3)}}(-u))$ in the left half plane. 
\end{corollary}

Observe that \eqref{pollaczekup} bypasses the cumbersome calculations with convolutions in \eqref{eq: thm2}. In \ref{ec:contour} we demonstrate that this is a numerically efficient way of computing the tight bounds. 

\subsection{Random walk lower bounds}


The tight upper bounds correspond to worst-case scenarios. We next show how the same MAD approach can identify best-case scenarios and hence tight lower bounds. For each $X_i$, define a second ambiguity set, which is a subset of 
$\mathcal{P}_{(\mu,d)}$:
\begin{equation}
\label{eq:mu_d_beta_conditions_multivariable}
\mathcal{P}_{(\mu,d,\beta)} = \left\{ \mathbb{P}: \ \mathbb{P} \in \mathcal{P}_{(\mu,d)},\ \mathbb{P} (X_i \geq \mu_i) = \beta_i, \ \forall i \right\}.
\end{equation}
Hence, for obtaining a lower bound, we include the additional information $\mathbb{P} (X_i \geq \mu_i) = \beta_i$ in the ambiguity set. Now, instead of finding the worst-case distribution, we want to identify the best-case distribution and corresponding tight lower bound. The following result is a direct consequence of the general lower bound in \cite{BenTal1972} on the expectation of a convex function of  independent random variables with $\mathcal{P}_{(\mu,d,\beta)}$ ambiguity. In Section~\ref{newproof3} we present a novel proof using the primal-dual method developed earlier for proving Theorem~\ref{thm1h}.


\begin{theorem}\label{thm1hj} 
\begin{equation}\label{lowww}
\min\limits_{\mathbb{P} \in \mathcal{P}_{(\mu,d,\beta)}} \expectp{h_n(\vec{X}) }= \sum\limits_{\bm{\alpha} \in \{1,2 \}^{n}} h_n(\upsilon_{\alpha_1}^{(1)} ,\ldots, \upsilon_{\alpha_{n}}^{(n)} )\prod\limits_{i=1}^{n} q_{\alpha_i}^{(i)}  ,
\end{equation}
where 
\begin{equation}
\label{eq:Best_case_relationships}
q_1^{(i)}  = \beta_i,\quad q_2^{(i)}  = 1-\beta_i,\quad \upsilon_1^{(i)}  = \mu_i + d_i/2\beta_i,\quad \upsilon_2^{(i)}  = \mu_i - d_i/2(1-\beta_i).
\end{equation}
\end{theorem}


Again specialize to the i.i.d.~setting, and denote by $Y$ the random variable with two-point distribution on values 
$$
v_1=\mu+\frac{d}{2\beta}, \quad v_2=\mu-\frac{d}{2(1-\beta)},
$$
with probabilities $\beta$ and $1-\beta$, respectively. Using similar reasonings as for the upper bound, we obtain for the tight lower bound for $\expect{M_n}$ an expression that sums over $O(n^2)$ terms:
\begin{equation}\label{permm12}
\sum_{k=1}^n \frac1k \sum_{k_1+k_2=k}\frac{k!}{k_1!k_2!}\beta^{k_1}(1-\beta)^{k_2}\max\{0,k_1v_1+k_2v_2\} .
\end{equation}
The tight lower bound for $c_m(M)$ can be expressed in terms of the integral
\begin{equation}\label{poll22}
\frac{(-1)^m}{2\pi i}\int_{\mathcal{C}}\frac{\log(1-\phi_{Y}(-u))}{u^{m+1}}{\rm d}u,
\end{equation}
where $\phi_Y(s)=\beta {\rm e}^{s v_1}+(1-\beta) {\rm e}^{s v_2}$,  $\mathcal{C}$ is a contour to the left of, and parallel to, the imaginary axis, and to the right of any singularities of $\log(1-\phi_{Y}(-u))$ in the left half plane. 

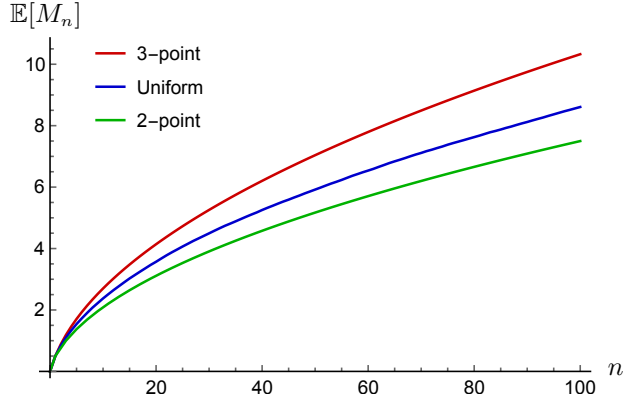
\begin{figure}
\begin{center}
\begin{tikzpicture}[scale=1.5]
\node[inner sep=0pt] (russell) at (0,0)
{\includegraphics[scale=0.6]{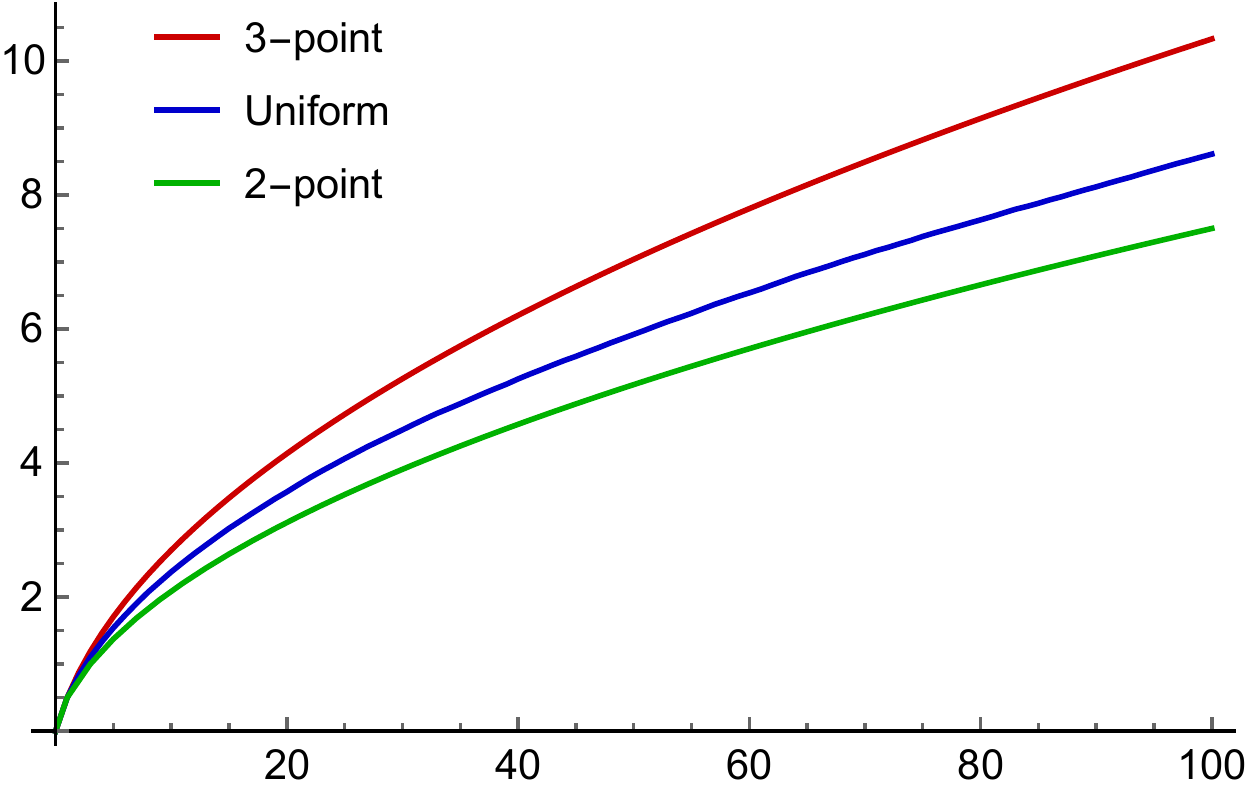}};
 \coordinate (a) at (2.7,-1.34);
 \coordinate (b) at (-2.35,1.8);
  \node at (a) {{\small $n$}};
  \node at (b) {{\small $\expect{M_n}$}};
\end{tikzpicture}
\caption{Expected random walk maximum $\expect{M_n}$ for $U(-b,b)$ and $b=2$ distributed step sizes with MAD $b/2$ (middle curve, obtained by simulation). The upper curve corresponds to the extremal three-point distribution within the ambiguity set with $\mu=0$, {$d=b/2$} and range $[-b,b]$, and the lower curve is the bound \eqref{permm12} from the two-point distribution with $\beta=1/2$.   
}\label{fig:EWu2}
\end{center}
\end{figure}

We illustrate the lower bound \eqref{lowww} (calculated using \eqref{permm12}) in Figure~\ref{fig:EWu2} for the random walk with step size $X$ having a uniform distribution on $[a,b]$. Here we assume a specific distribution just for illustration purposes. The MAD of $X$ can be shown to be $(b-a)/4$. In Figure~\ref{fig:EWu2} we choose $b=-a=2$ so that $\mu=0$ and $d=1$. Observe that upper and lower bound together  provide a tight interval for all possible distributions in the ambiguity set $\mathcal{P}_{(0,1,1/2)}$.

\begin{figure}
\begin{center}
\begin{tikzpicture}[scale=1.5]
\node[inner sep=0pt] (russell) at (0,0)
{\includegraphics[scale=0.6]{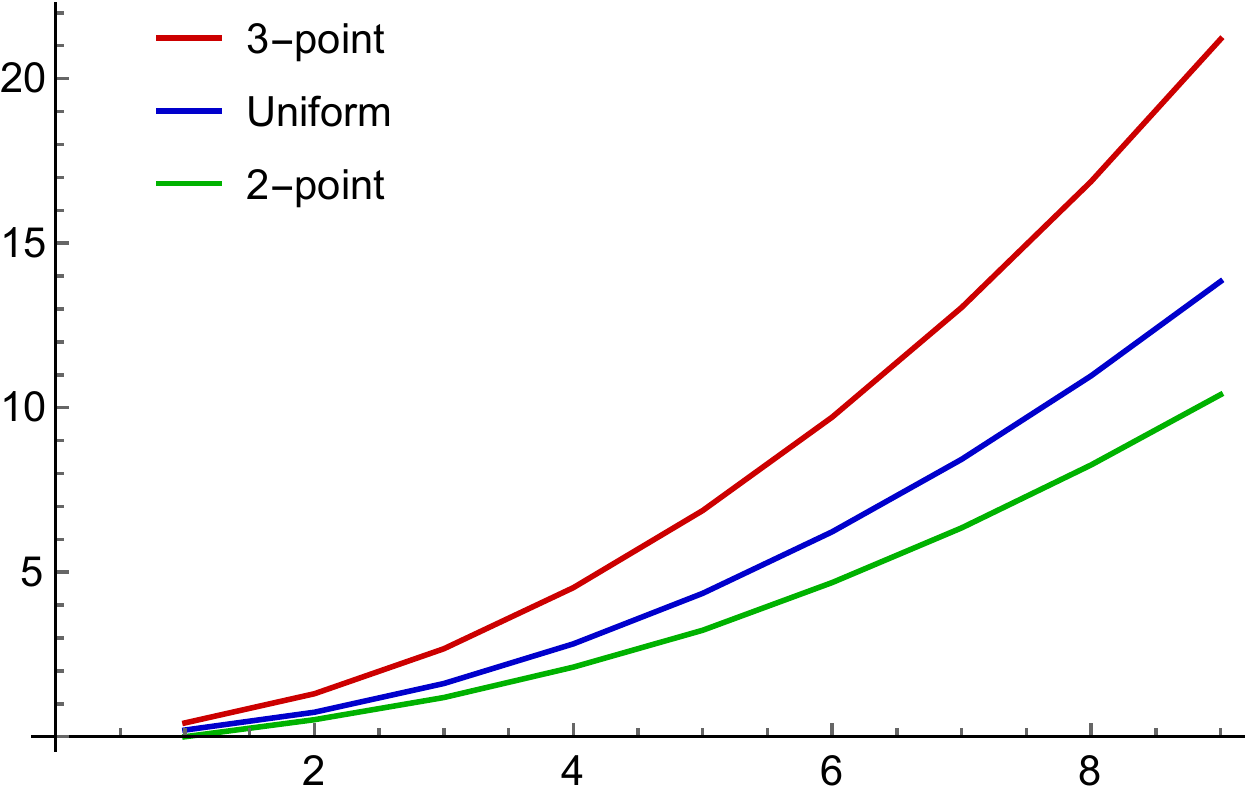}};
 \coordinate (a) at (2.7,-1.34);
  \coordinate (b) at (-2.35,1.8);
  \node at (a) {{\small$b$}};
    \node at (b) {{\small $\expect{M}$}};
\end{tikzpicture}
\caption{Expected all-time maximum $\expect{M}$ for $U(-b-2,b)$ and $b\in(1,10)$ (middle curve, obtained by simulation). The upper curve corresponds to the extremal three-point distribution within the ambiguity set with $\mu=-1$,  {$d=(b+1)/2$} and range $[-b-2,b]$, and the lower curve is the bound \eqref{permm12} from the two-point distribution with $\beta=1/2$.   
}\label{fig:EWu}
\end{center}
\end{figure}
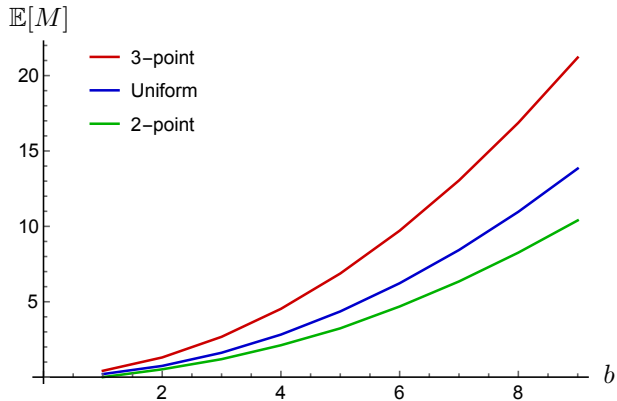

Figure~\ref{fig:EWu} shows the tight upper bound \eqref{pollaczekup} and the lower bound \eqref{poll22} for 
 $\expect{W}$ with ambiguity set with $\mu=-1$, $d=b/2$ and range $[-b-2,b]$. Observe that the bounds increase with the range and the MAD (which can be shown to hold in general). For a point of reference, we also plot the exact results for one member of the ambiguity set, when generic increment having a uniform distribution on $[-b-2,b]$. 


{\subsection{Comparison with mean-variance ambiguity}\label{secmvc}
As explained earlier, mean-variance ambiguity appears less computationally tractable than mean-MAD ambiguity. We now show how the key result for mean-MAD ambiguity, Theorem~\ref{thm1h}, can be used to obtain results for mean-variance ambiguity. Let $\mathcal{P}^*_{(\mu,\sigma)}$ denote the ambiguity set that contains all distributions with known range, mean and variance, i.e.
\begin{equation}
\label{}
\mathcal{P}^*_{(\mu,\sigma)} = \left\{ \mathbb{P}: \ \text{supp}(X_i) \subseteq [a_i,b_i], \ \mathbb{E}_{\mathbb{P}}(X_i) = \mu, \ \mathbb{E}_{\mathbb{P}}( X_i - \mu )^2 = \sigma^2,\ \forall i, \ X_i \indep X_j, \ \forall i \neq j \right\}. 
\end{equation}

\begin{proposition}\label{propp}
Let $d_{\rm min}=2\sigma^2/(b-a)$ and $d_{\rm max}=\sigma$. Then,
\begin{equation}\label{cvb}
\max\limits_{\mathbb{P} \in \mathcal{P}_{(\mu,d_{\rm min})}} \expectp{h_n(\vec{X})}\leq \max\limits_{\mathbb{P} \in \mathcal{P}^*_{(\mu,\sigma)}} \expectp{h_n(\vec{X})}\leq \max\limits_{\mathbb{P} \in \mathcal{P}_{(\mu,d_{\rm max})}} \expectp{h_n(\vec{X})}
\end{equation}
\end{proposition}
\proof{Proof.}
From \cite{BenTal1985}, we know that
$$
\frac{2 \sigma^2}{b-a}\leq d\leq \sigma. 
$$
Hence, $
\max_{\mathbb{P} \in \mathcal{P}^*_{(\mu,\sigma)}} \expectp{h_n(\vec{X})}=
\max_{\mathbb{P} \in \mathcal{P}_{(\mu,d^*)}} \expectp{h_n(\vec{X})}$
for some $d^*\in [2\sigma^2/(b-a),\sigma]$. Since $\max_{\mathbb{P} \in \mathcal{P}_{(\mu,d)}} \expectp{h_n(\vec{X})}$ is non-decreasing in $d$, see \cite{postek2018robust}, the result follows.  
\Halmos
\endproof

Notice that Proposition~\ref{propp} presents a way to delimit the upper bounds of all stationary cumulants $c_m(M)$ and the transient mean $\expect{M_n}$ under mean-variance ambiguity. The mean-MAD bounds are specified in terms of specific three-point distributions. 

We next show that the lower bound  in Proposition~\ref{propp} can lead to a result for infinite-support distributions. Set $b=a+\xi(\mu-a)$ with $\xi\geq 1$, and observe that the lower bound in \eqref{cvb} comes with the extremal three-point distribution
$$
X_{(3)}^\xi=\left\{\begin{array}{ll}
a &  {\rm w.p.} \  \frac{\sigma^2}{(\mu-a)^2\xi},\\
\mu &  {\rm w.p.} \  1-\frac{\sigma^2}{(\mu-a)^2\xi}-\frac{\sigma^2}{(\mu-a)^2\xi(\xi-1)},\\
a+\xi(\mu-a) &  {\rm w.p.} \  \frac{\sigma^2}{(\mu-a)^2\xi(\xi-1)}.\\
  \end{array}\right.
$$
This distribution has mean $\mu$ and variance $\sigma^2$, irrespective of the range $[a,b]$. We can thus let $\xi$ grow to infinity to investigate what happens for infinite-support distributions.

For the expected all-time maximum, we can exploit an argument very similar to \cite{chen2019extremal}, Theorem EC.3. A classic result from regenerative analysis says that the expected all-time maximum is the expected sum of the random walk position over one cycle, denoted by $\E[{\rm integral}]$, divided by the expected length of one cycle, i.e. $\E[{\rm cycle \ length}]$. This cycle will consists of a period during which the queue remains empty, corresponding of consecutive (negative) steps of size $a$ or $\mu$.
As $\xi$ increases, the three-point distribution places probabilities of order $O(1/\xi^2)$ on $a$ and $a+\xi(\mu-a)$, and the rest of the mass on point $\mu$. 
As $\xi$ grows large, only rarely with probability $O(1/\xi^2)$, a large positive step occurs. The impact of the very large step of size $a+\xi(\mu-a)$ is 
roughly the area of the triangle with height $a+\xi(\mu-a)$ and width $(a+\xi(\mu-a))/(-\mu)$, and hence $\E[{\rm integral}]=(a+\xi(\mu-a))^2/(-2\mu)\sim(\xi(\mu-a))^2/(-2\mu)$ as $\xi \to\infty$. The cycle then consists of an empty period of expected length $(1-p_b)/p_b\sim (\xi(\mu-a))^2/\sigma^2$ and the positive period due to the large step of  expected length $(a+\xi(\mu-a))/(-\mu)$, so that  $\E[{\rm cycle \ length}]\sim (\xi(\mu-a))^2/\sigma^2$, and the expected all-time maximum converges to $\sigma^2/(-2\mu)$ as $\xi\to\infty$. Since this is a lower bound for $\max_{\mathbb{P} \in \mathcal{P}^*_{(\mu,\sigma)}}\expect{M}$, we know that for the random walk with generic step size $X$ it holds that 
$\max_{\mathbb{P} \in \mathcal{P}^*_{(\mu,\sigma)}}\expect{M}\geq \sigma^2/(-2\mu)$. This lower bound matches Kingman's upper bound  $\expect{M}\leq \sigma^2/(-2\mu)$, which proves that Kingman's upper bound is tight. 
Tightness of Kingman's bound was already proven in \cite{daley1992inequalities} by identifying a two-point distribution with mean $\mu$, variance $\sigma^2$ such that  $\expect{M}$ approaches the upper limit as one of the two points goes to infinity. 





}

{
\subsection{Degenerate behavior for infinite range}

Compared to variance, MAD may be more appropriate in case of real-life empirical data that display non-Gaussian features and outliers. Indeed, unlike standard deviation, MAD does not require existence of second moments, and is not so much affected by large deviations from the mean. This feature, however, has major consequences when we let the range $[a,b]$ grow large in which case conditioning on the MAD being $d$ thus allows for distributions with relatively heavy tails. In particular, in the limit $b\to\infty$, this will lead to overly pessimistic scenarios as heavy-tailed distributions with infinite second moments would still have a finite $d$ and hence be member of the ambiguity set. While for large but finite $b$ a truly heavy-tailed distribution with infinite second moment is ruled out, the dispersion allowed by the ambiguity set might become too loose for practical purposes.  
 An effective usage of the robust mean-MAD framework therefore requires a careful selection of the range, for which we now present some guidelines. 

Observe that the variance of $X_{(3)}$ is $\frac{d}{2}(b-a)$, the maximal variance for distributions in the ambiguity set $\mathcal{P}_{(\mu, d)}$. Hence, for fixed $d$, the variance becomes unbounded when $b\to\infty$. As a consequence, this results in fairly crude bounds:
\begin{proposition}
As $b\to\infty$, the bound  $\max\limits_{\mathbb{P} \in \mathcal{P}_{(\mu,d)}} \expectp{f_n(\vec{X})}$ converges to
\begin{equation}\label{bounddon}
n\cdot\frac{d}{2} + \sum_{k=1}^n \frac{1}{k} \sum_{k_1+k_2=k}\max\{0,k_1a+k_2\mu\}\cdot \frac{k!}{k_1!k_2!}p_1^{k_1}p_2^{k_2}
\end{equation}
with $p_1 = \frac{d}{2(\mu-a)}$ and  $p_2 = 1 - \frac{d}{2(\mu-a)}$.
\end{proposition}
\proof{Proof.}
Split the inner 
summation in \eqref{permm6} into three parts. 
First consider the summation over $\sum_i k_i=k:k_3\geq2$, hence those instances for which the value $b$ occurs multiple times. Taking the limit $b\to\infty$ inside of the summation and recognizing the fact that the probability mass on the third point is $O(\frac{1}{b^{k_3}})$ gives
    \begin{equation}
        \lim\limits_{b\to\infty} \sum_{\sum_i k_i=k:k_3\geq2}\frac{d^{k_3}\max\{0,k_1a+k_2\mu+k_3b\}}{2^{k_3}(b-\mu)^{k_3}}\cdot \frac{k!}{k_1!k_2!k_3!}p_1^{k_1}p_2^{k_2} = 0.
    \end{equation}
    Next consider ${\sum_i k_i=k:k_3=1}$, describing the instances for which the extremal point $b$ occurs precisely once. Taking the limit $b\to\infty$ inside the sum and using that the probability mass on the point $b$ is $O(\frac{1}{b})$ gives
    \begin{equation}
        \lim\limits_{b\to\infty} \sum_{\sum_i k_i=k:k_3=1}\frac{d\max\{0,k_1a+k_2\mu+b\}}{2(b-\mu)}\cdot \frac{k!}{k_1!k_2!}p_1^{k_1}p_2^{k_2} = k\cdot\frac{d}{2} \cdot \sum_{\sum_i k_1+k_2=k-1} \frac{(k-1)!}{k_1!k_2!}p_1^{k_1}p_2^{k_2} = k\cdot\frac{d}{2}.
    \end{equation}
    The third part is then ${\sum_i k_i=k:k_3=0}$, representing the instances without occurrence of the point $b$. Taking the limit inside of the summation we get
    \begin{equation}
        \lim\limits_{b\to\infty} \sum_{\sum_i k_i=k:k_3=0}\max\{0,k_1a+k_2\mu\}\cdot \frac{k!}{k_1!k_2!}p_1^{k_1}p_2^{k_2} = \sum_{k_1 + k_2=k}\max\{0,k_1a+k_2\mu\}\cdot \frac{k!}{k_1!k_2!}p_1^{k_1}p_2^{k_2}.
    \end{equation}
    This completes the proof. 
\Halmos
\endproof

The proof reflects that large running maxima are likely due to a single large step. The feature is caused by heavy-tailed distributions, and in queueing theory dubbed the single big jump principle (see e.g., \cite{Foss2007}). This dominance of one step sharply contrasts intuition for light-tailed distributions, where typically all steps together lead to large sums or maxima. 
The bound \eqref{bounddon} for $\expect{M_n}$ grows to infinity as $n\to\infty$, rendering the bound useless for the expected all-time maximum $\expect{M}$. This is indeed anticipated, and can be understood as follows. Define a sequence of random walks indexed by $b$ with the extremal three-point distribution. Consider the limiting all-time maximum $M$ as $b\to \infty$. Assume that the random walk has negative drift (i.e., $\expect{X}<0$). Then the associated sequence of distributions of $M=M_{(b)}$ will converge to a proper limit $M_{(\infty)}$. However, as $\lim_{b\to\infty}\mathcal{P}_{(\mu,d)}$ contains distributions with infinite second moment, \cite{asmussen2003}, Theorem~\MakeUppercase{\romannumeral10}.2.1, says that $\expect{M_{(\infty)}}$ will be infinite.


\subsection{Setting the range to construct adequate bounds}

 We now present some guidelines for setting the range, based on the observation that many distributions come with a MAD and standard deviation of comparable size. For the Pearson family of distributions (which includes the gamma and normal distribution) with mean $\mu$ and variance $\sigma^2$, the MAD $d$ and variance are related as
\eqan{
d=2 \alpha \sigma^2 p(\mu)
}
with $\alpha$ a constant depending on skewness and kurtosis and $p(\mu)$ the density in $\mu$. For the exponential distribution this relation gives $d=(2/\e)\sigma$ and for the normal distribution $d=(\sqrt{2/\pi})\sigma$. Other distributions for which the ratio $d/\sigma$ is constant include the uniform distribution and discrete distributions such as the Poisson, binomial, and negative binomial distribution. 
With this in mind, in a way similar to constructing confidence intervals in statistical estimation, we then choose to set the range as the mean plus or minus a constant times the MAD:
\begin{equation}\label{eq:truncation}
    a = \mu - k\cdot d, \quad b = \mu + k\cdot d.
\end{equation}
Here we regard $d$ as the natural scale of deviation, and $k$ as a free parameter that sets the robustness level. So we take the mean and MAD as given, and regard the range as tunable  (using common sense or statistical evidence)  by the decision maker. We should stress that, while intuitive from a probabilistic perspective,  the rule \eqref{eq:truncation} is only one of many ways to choose the parameters $a,b$. 


We demonstrate \eqref{eq:truncation} for a setting where we take the M/M/1 queue as the `true' model. The increment $X$ now becomes the difference of two exponential random variables for which we have a closed-form MAD expression in terms of the mean value of $X$ (see the caption of Table~\ref{res:trunc}). We thus have reference values for $\mu$ and $d$, and can investigate the impact of $k$.
Observe that the bound grows almost linearly with $k$, in particular in heavy-traffic scenarios, and this underlines the need for careful selection of the range. While the actual range of the M/M/1 queue spans all real numbers, we see that restricting deviations to twice the MAD ($k=2$) gives comparable model performance. 
When reading Table~\ref{res:trunc}, keep in mind that the overall goal in this paper is not to approximate specific models (such as the M/M/1 queue), but rather to come with conservative, robust estimates for an entire class of models that share the same mean-MAD-range properties. In that sense, $k=2$ is not better than $k=1.5$ or $k=2.5$, but rather expresses a different ambiguity assessment or robustness level. 

\begin{table}[h!]
\caption{The actual value and bounds of the expected steady-state waiting time
$\expect{W}$ of the M/M/1 queue with unit mean exponential interarrival times and exponential service times with mean $\rho$, where the increment $X$ has mean $\mu=\rho-1$ and MAD $d=\frac{2e^{\rho-1}}{\rho+1}$, with the range $[a,b]$ set through the rule \eqref{eq:truncation}.}\label{res:trunc}
 \vspace{.3cm}
\begin{center}
\begin{tabular}{rrrrrrrr}\hline
     &          & \multicolumn{6}{c}{$k$} \\ \cline{3-8}
$\rho$     &   $\mathbb{E}[W]$       & 1.5      & 1.75     & 2  & 2.25     & 2.5       & 3         \\ \hline
0.1  & 0.01111  & 0.10497  & 0.16434  & 0.21535  & 0.25915   & 0.30116   & 0.40782   \\
0.5  & 0.50000  & 0.56329  & 0.67919  & 0.79663  & 0.91459   & 1.02840   & 1.26462   \\
0.6  & 0.90000  & 0.86690  & 1.03323  & 1.19770  & 1.36332   & 1.52804   & 1.85818   \\
0.7  & 1.63333  & 1.41436  & 1.66589  & 1.91885  & 2.17142   & 2.42373   & 2.92850   \\
0.8  & 3.20000  & 2.57273  & 3.01339  & 3.45454  & 3.89573   & 4.33672   & 5.21866   \\
0.9  & 8.10000  & 6.21057  & 7.25250  & 8.29428  & 9.33642   & 10.37811  & 12.46184  \\
0.99 & 98.01000 & 73.55537 & 85.81540 & 98.07540 & 110.33542 & 122.59543 & 147.11548 \\ \hline
\end{tabular}
\end{center}
\end{table}

}

\section{Extremal GI/G/1 queue}\label{sec4}
Let us now turn to the extremal GI/G/1 queue problem, as described in the introduction. 
Let $W_n$ be the waiting time of customer $n$. The sequence $(W_n, n \geq 0)$ with $W_0=0$ satisfies the Lindley recursion
\begin{equation}\label{thelind}
W_{n+1} =(W_n +{\VV}_n-{\UU}_n)^+, \quad  n\geq 0.
\end{equation}
Let $W$ be the steady-state waiting time. Since $W_n \das M_n$ and $W \das M$ the results for the random walk maxima likely carry over to the waiting times. The main difference is that the step size  $X$ is now interpreted as the difference $\VV-\UU$ between the generic service time and generic interarrival time. If one has mean-MAD information about both $\VV$ and $\UU$ this is more informative than mean-MAD information about $\VV-\UU$, and this additional information should lead to even sharper bounds. 

\subsection{A complete picture}

The GI/G/1 queue assumes that interarrival times and service times are independent, so it is natural to assume that 
$\VV$ has ambiguity set $\mathcal{P}_{(\mu_{\VV}, d_{\VV})}$ and  $\UU$ has ambiguity set $\mathcal{P}_{(\mu_{\UU}, d_{\UU})}$, where the ambiguity sets now contain all distributions for \textit{univariate} $\VV$ and $\UU$, that is,
$$
\mathcal{P}_{(\mu_{\VV},d_{\VV})} = \left\{ \mathbb{P}: \ \text{supp}(\VV) \subseteq [a_{\VV},b_{\VV}], \ \mathbb{E}_{\mathbb{P}}(\VV) = \mu_{\VV}, \ \mathbb{E}_{\mathbb{P}}\left| \VV - \mu_{\VV} \right| = d_{\VV}\right\}
$$
and
$$
\mathcal{P}_{(\mu_{\UU},d_{\UU})} = \left\{ \mathbb{P}: \ \text{supp}(\UU) \subseteq [a_{\UU},b_{\UU}], \ \mathbb{E}_{\mathbb{P}}(\UU) = \mu_{\UU}, \ \mathbb{E}_{\mathbb{P}}\left| \UU - \mu_{\UU} \right| = d_{\UU}\right\}.
$$
The extremal queue problem with mean-MAD dispersion information can then be phrased as
\begin{equation}\label{maxphrase2}
\max\limits_{\mathbb{P} \in \mathcal{P}_{(\mu_{\VV},d_{\VV})}\times\mathcal{P}_{(\mu_{\UU},d_{\UU})}} \expect{ f (\vec{X})},
\end{equation}
where $\expect{ f (\vec{X})}$ describes $\expect{W_n}$ or  $c_m(W)$ and $\vec{X}$ is the random vector with elements $\UU_1,\VV_1,\UU_2,\VV_2,\ldots$. 
This is the classical setting of the extremal GI/G/1 queue treated in \cite{rolski1972some,eckberg1977sharp,whitt1984approximations,
chen2019extremal}, but with MADs instead of variances describing the ambiguity set. 
Let the random variables $\VV_{(3)}$ and $\UU_{(3)}$  follow the extremal three-point distributions $\Omega(\mu_{\VV},d_{\VV},a_{\VV},b_{\VV})$ and  $\Omega(\mu_{\UU},d_{\UU},a_{\UU},b_{\UU})$, respectively.
\begin{theorem}\label{tightgg1}
Consider the GI/G/1 queue with generic interarrival time 
$\UU$ with ambiguity set $\mathcal{P}_{(\mu_{\UU}, d_{\UU})}$ and generic service times  $\VV$ with ambiguity set $\mathcal{P}_{(\mu_{\VV}, d_{\VV})}$. Consider the tight upper bounds for the transient mean waiting time $\expect{W_n}$ and all cumulants of the steady-state waiting time $W$. 
\begin{hitemize}
\item[{\rm(i)}] For given interarrival time $\UU$, the tight upper bounds follow from the  service time $\VV_{(3)}$.
\item[{\rm(ii)}] For given service time $\VV$, the tight upper bounds follow from  the interarrival time $\UU_{(3)}$. 
\item[{\rm (iii)}] The overall tight upper bounds follow from interarrival time $\UU_{(3)}$ and service time $\VV_{(3)}$.
\end{hitemize}
\end{theorem}
\proof{Proof.}
Like Theorem~\ref{thm1h}, the tight bounds for $\expect{W_n}$
follow from the general upper bound in \cite{BenTal1972} on the expectation of a convex function of the random vector $(X_1,\ldots,X_n)$ with mean-MAD ambiguity, but now with $X_i$  replaced by $\VV_i-\UU_i$. The function describing $\expect{W_n}$ (see Theorem~\ref{thm1h}) is indeed convex in {\it both} $\VV_i$ and $\UU_i$, and hence the result follows. Similarly, Spitzer's formula for $c_m(W)$  (see Theorem~\ref{thm1f})  is also convex in both $\VV_i$ and $\UU_i$, and hence the tight bounds for $c_m(W)$ follow from our proof of Theorem~\ref{thm1f}.
 \Halmos
\endproof

Using the earlier results for the random walk, we present in \ref{ec:contour} expressions that are  helpful in evaluating the tight bounds.  
Table \ref{tab11}
shows an example of the tight bound for $\expect{W}$ associated with $(\UU_{(3)}, \VV_{(3)})$, also compared with other known bounds that require variance information (see \ref{secdfb}). 
The variance of the extremal three-point distribution $\Omega(\mu,d,a,b)$ is $\frac{d}{2}(b-a)$, the maximal variance for distributions in the ambiguity set $\mathcal{P}_{(\mu, d)}$. We thus know the variances of $\UU_{(3)}$ and $\VV_{(3)}$, and can calculate the other three bounds. 
In heavy traffic, Kingman's bound is known to be asymptotically correct, and hence the other three (sharper) bounds also converge to the heavy-traffic limit as $\rho\uparrow 1$. See \ref{sec:ECn} for more numerical results. {Notice that Table \ref{tab11} is not meant to compare mean-MAD with mean-variance bounds. The displayed differences merely express different ways of dealing with ambiguity. Also remember that the mean-MAD bounds in Theorem~\ref{tightgg1} are crucially influenced by the choice of range, in this example set to $[0,10]$ for both the interarrival and service time distributions.}

\begin{table}[h]
\caption{{Bounds for $(1-\rho)\expect{W}/\rho$  for $(\mu_{\UU},d_{\UU},a_{\UU},b_{\UU})=(1,1,0,10)$ and 
 $(\mu_{\VV},d_{\VV},a_{\VV},b_{\VV})=(\rho,0.1,0,10)$.} }\label{tab11}
 \vspace{.3cm}
\begin{center}
\begin{tabular}{r C{1.5cm} R{1.5cm} R{1.5cm} R{1.5cm}}\hline
   $\rho$    & Thm.~\ref{tightgg1} & C \& W \eqref{chenwhittbound} & Daley \eqref{daley} & Kingman \eqref{kingman}  \\ \hline
0.1  & 4.06613 & 7.00020 & 7.25000 & 27.50000 \\
0.2  & 2.52306 & 5.27810 & 5.75000 & 13.75000 \\
0.5  & 2.03141 & 3.63750 & 4.25000 & 5.50000  \\
0.7  & 2.49160 & 3.17138 & 3.60714 & 3.92857  \\
0.8  & 2.61932 & 3.00523 & 3.31250 & 3.43750  \\
0.9  & 2.69802 & 2.86711 & 3.02778 & 3.05556  \\
0.95 & 2.72609 & 2.80627 & 2.88816 & 2.89474  \\
0.99 & 2.74547 & 2.76091 & 2.77753 & 2.77778  \\ \hline
\end{tabular}
\end{center}
\end{table}






\subsection{Further comparison between MAD and variance}
For the variance counterpart, \cite{chen2019extremal} also formulate a semi-infinite linear optimization problem. The crucial difference is that they cannot use the univariate function extension (as explained in Section~\ref{sec2}), and hence should work directly with the multivariate function. This in turn implies that the dual problem cannot be solved explicitly (like in the univariate case), let alone that there is a zero duality gap. 
Another complication is that the multivariate function based on Spitzer's formulas \eqref{permm3} and \eqref{permm6} cannot be expressed directly in $\VV$ and $\UU$, but rather in terms of convolutions of the distributions of $\VV$ and $\UU$. 
\cite{chen2019extremal} resolve these considerable challenges by several ingenious arguments, a.o.~exploiting the description of $W$ as a fixed point of the stochastic equation $W\das (W+\VV-\UU)^+$, and by imposing additional regularity conditions on $\VV$. In this way,  \cite{chen2019extremal} prove a similar but weaker result than 
 Theorem~\ref{tightgg1} for the exact same setting, but with variance as dispersion measure. They show that the extremal distributions of $\UU$ and $\VV$ both have supports on at most three points. 

An important message of this paper is that with MAD the extremal distribution remains unaltered going from the univariate to the multivariate setting, and that with variance this reasoning fails. In fact, one intuitively expects formidable challenges when seeking for  extremal distributions under variance constraints. This intuition is confirmed by Chen and Whitt's formulation of the extremal distribution as the solution of a non-convex nonlinear optimization problem. While this optimization problem can be solved numerically, 
a closed-form solution and hence identification of the extremal distribution remains out of reach.  

Under variance constraints, it is conjectured that the tight bound comes from specific two-point distributions for both $\UU$ and $\VV$. In fact, the bound \eqref{chenwhittbound} in Table~\ref{tab11} holds under the assumption that this conjecture is true, and was shown by \cite{chen2019extremal}  to be very close to the tight upper bound. Theorem \ref{tightgg1} rules out a similar conjecture in the MAD setting. The tight bounds in  Theorem \ref{tightgg1}  always involve three-point distributions.



{
\subsection{Data-driven setting}

In applications, you may only have a limited number $n$ of observed interarrival and service times. We consider this realistic setting where knowledge of the stochastic nature is
restricted to a set of samples generated independently and randomly according to an unknown distribution $\P$. To apply the mean-MAD framework in this context, we need to construct the ambiguity set that is supposed to contain this unknown $\P$. We will show that we can efficiently estimate the mean, MAD, and $\beta$, and hence compute robust bounds that are useful in realistic settings.

Let $\mu_n^{(V)}$, $d_n^{(V)}$ and $\beta_n^{(V)}$ denote the consistent estimators of $\mu_V$, $d_V$ and $\beta_V=\prob{V\geq\mu_V}$, respectively, based on $n$ observed service times $v_1,...,v_n$, and defined as $\mu^{(V)}_n=\bar{v}=\frac1n \sum_{i=1}^n v_i$,  $d^{(V)}_n=\frac1n\sum_{i=1}^n |v_i-\bar{v}|$ and  $\beta^{(V)}_n=\frac1n\sum_{i=1}^n \mathbbm{1}_{[\bar{v},\infty)}(v_i)$ as consistent estimators.
We define similar estimators based on $n$ observed interarrival times. 
Next, we demonstrate the mean-MAD bounds in this data-driven setting.  Since statistical accuracy of the estimators increases with the number of samples, we expect the bounds to converge as $n$ increases. Figure \ref{fig:dataconvergence} illustrates two sample paths representing the estimates for the upper and lower bound and their convergence to the tight mean-MAD bounds, where $V$ and $U$ both follow a uniform distribution on the intervals $[0,5]$ and $[0,10]$, respectively. Observe that convergence settles in quickly.



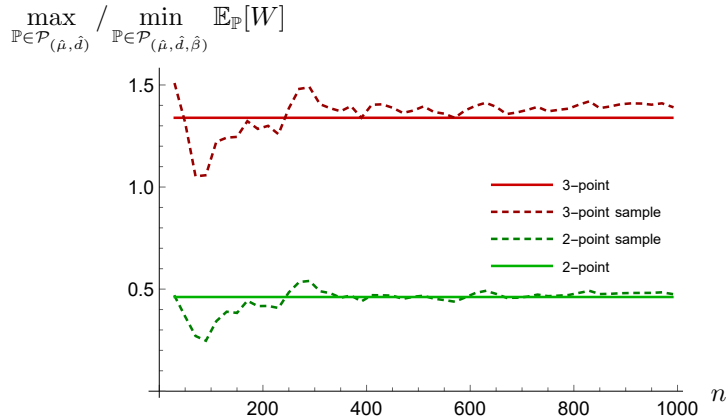
\begin{figure}[h!]
\begin{center}
\begin{tikzpicture}[scale=1.5]
\node[inner sep=0pt] (russell) at (0,0)
{\includegraphics[scale=0.6]{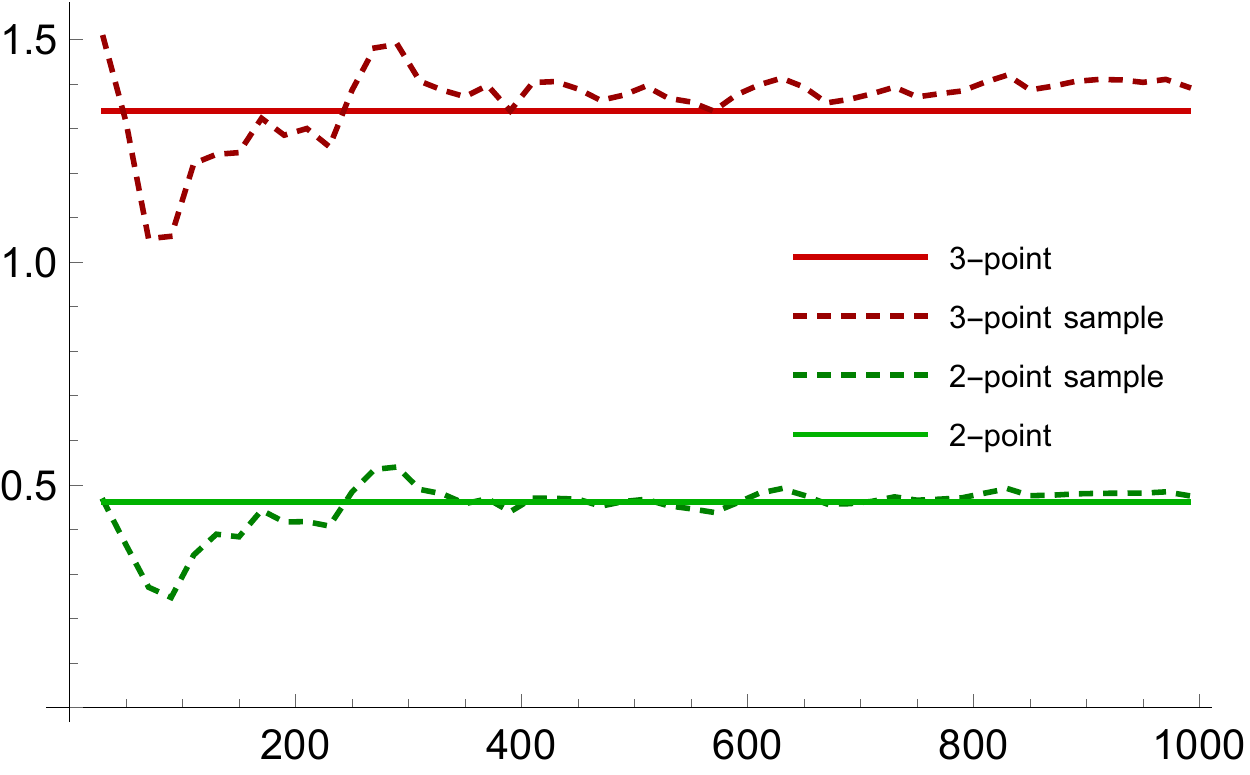}};
 \coordinate (a) at (2.7,-1.34);
 \coordinate (b) at (-2.35,1.9);
  \node at (a) {{\small $n$}};
  \node at (b) {{\small $\max\limits_{\mathbb{P}\in\mathcal{P}_{(\hat{\mu},\hat{d})}}/\min\limits_{\mathbb{P}\in\mathcal{P}_{(\hat{\mu},\hat{d},\hat{\beta})}}\expectp{W}$}};
\end{tikzpicture}
    \caption{Estimation of the mean-MAD ambiguity upper and lower bound. The red and green line represent the true upper and lower bound, respectively, and the dashed lines represent bound estimates which are computed using the realizations $v_1,...,v_n$ drawn from a $U(0,5)$ distribution and $u_1,...,u_n$ sampled from a $U(0,10)$ distribution.}
    \label{fig:dataconvergence}
\end{center}
\end{figure}



We have also performed extensive simulations to investigate the error between the estimated and true bounds for several values of the sample size $n$. We generate 1,000 sample paths of sample size 10,000 and compute the corresponding mean relative error. Table \ref{res:dataexperiment} displays the mean absolute percentage error (MAPE) for both the upper and lower bound estimates, where the interarrival time is $U(0,10)$ distributed and we differentiate between a 50\% and 90\% utilization level. Observe that estimating the lower bound is slightly harder than estimating the upper bound. Indeed, the lower bound requires estimating the additional parameters $\beta_V$ and $\beta_U$. 
Also observe that the relative error increases with the system utilization.

\begin{table}[h!]
\caption{MAPE of the bound estimates for $n\in\{150, 200, 500, 1000, 2000, 5000, 10000\}$. The interarrival times are $U(0,10)$ distributed and the results differentiate between two service time distributions and the upper and lower mean-MAD bounds. Sample paths resulting in instable systems were removed and done over.}\label{res:dataexperiment}
 \vspace{.3cm}
\begin{center}
\begin{tabular}{ccccccccc}\hline
     &          &  \multicolumn{7}{c}{MAPE with sample size $n$} \\ \cline{3-9}
  Service times    &  Bound  & 150      & 200     & 500     & 1000     & 2000    & 5000 & 10000   \\ \hline
$U(0,5)$ & UB & 15.44\% & 13.22\% & 8.31\%  & 5.84\%  & 4.28\%  & 2.72\%  & 1.89\% \\
      & LB & 25.51\% & 22.30\% & 13.86\% & 9.75\%  & 7.08\%  & 4.53\%  & 3.16\% \\
$U(0,9)$ & UB & 33.35\% & 30.93\% & 21.93\% & 16.35\% & 13.29\% & 8.92\%  & 6.41\% \\
      & LB & 36.27\% & 35.01\% & 28.72\% & 22.30\% & 17.35\% & 10.77\% & 7.58\%\\ \hline
\end{tabular}
\end{center}
\end{table}

To further highlight the role of system utilization, we perform a similar data-driven experiment, but now with ground truth a single trace of $n$ customers in an  M/M/1 queue. The results are shown in Figure~\ref{fig:dataconvergenceLB}. Indeed, as $\rho$ increases, more observations are required for accurate parameter estimates and hence accurate bounds.

\begin{figure}[h!]
\begin{center}
\begin{tikzpicture}[scale=1.5]
\node[inner sep=0pt] (russell) at (0,0)
{\includegraphics[scale=0.6]{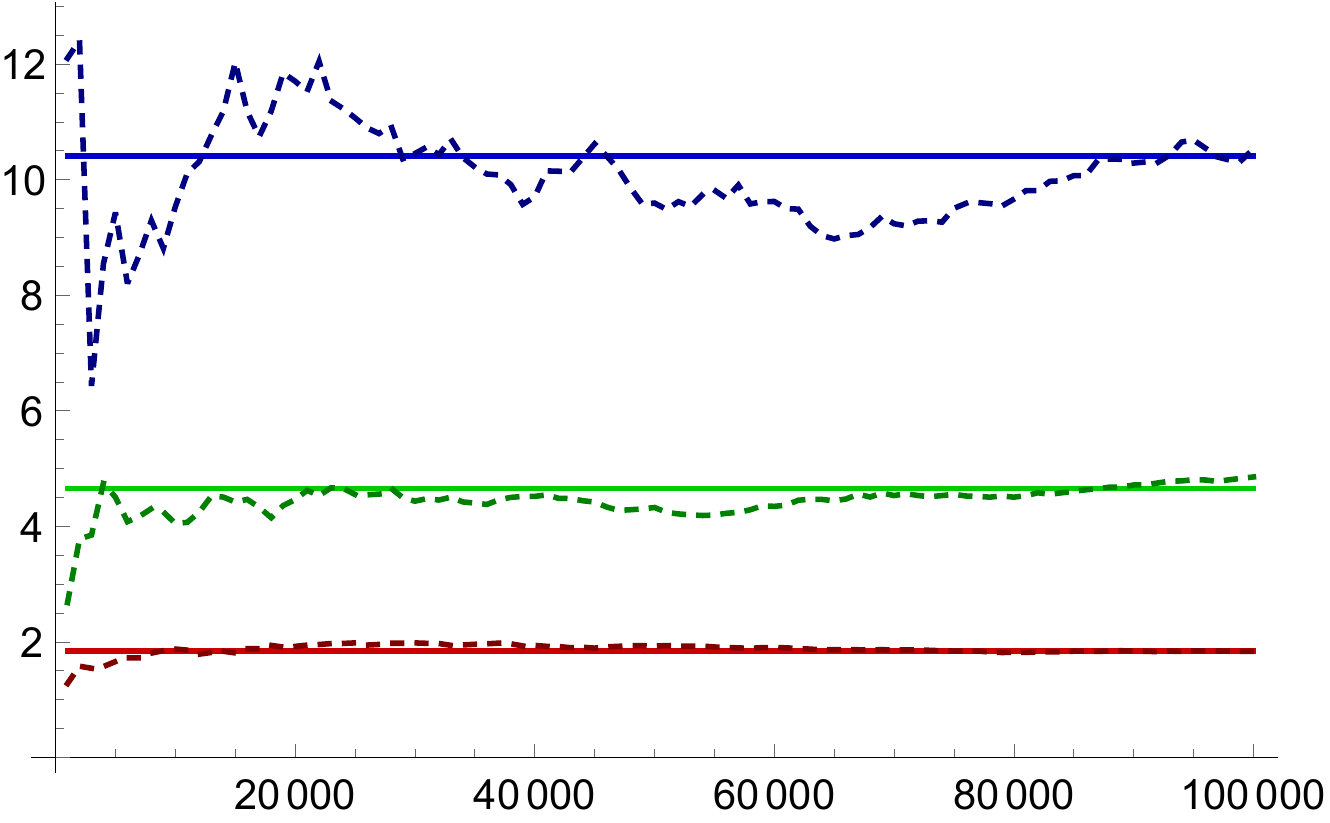}};
 \coordinate (a) at (2.7,-1.34);
 \coordinate (b) at (-2.35,1.9);
  \node at (a) {{\small $n$}};
  \node at (b) {{\small $\min\limits_{\mathbb{P}\in\mathcal{P}_{(\hat{\mu},\hat{d},\hat{\beta})}}\expectp{W}$}};
\end{tikzpicture}
    \caption{Estimation of the mean-MAD ambiguity lower bound for the M/M/1 queue. The solid red, green, and blue lines depict the bounds for $\rho=0.8,0.9,0.95$, respectively. The dashed lines represent the corresponding estimates of the bounds, where the $U_i$ are sampled from a unit mean exponential distribution and the $V_i$ are exponentially distributed with mean $\rho$.}
    \label{fig:dataconvergenceLB}
\end{center}
\end{figure}
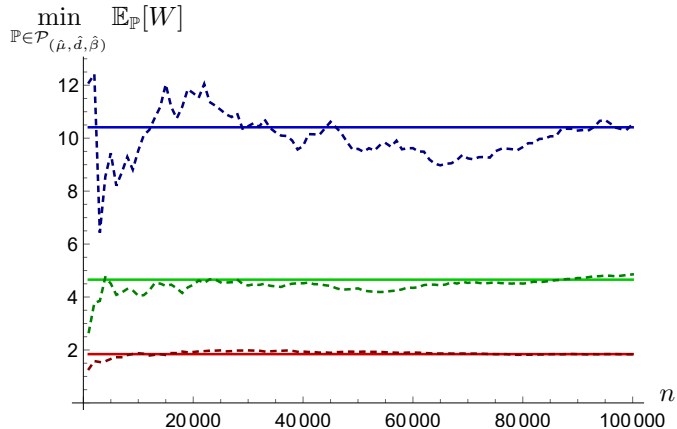
 

Taken together, we conclude that the robust bounds are useful for realistic data-driven settings that require statistical estimation of the summary statistics such as the mean and MAD.

}

\section{Conclusions}\label{sec6}


This paper explains why MAD simplifies comparable variance-based optimization problems, in a way that is almost unreasonably effective, resulting in  a full solution to the extremal queue problem with mean-MAD constraints. 
When partial information is available in the form of mean, range and MAD, we have obtained the sharpest possible  bounds. Through basic statistical estimation of this partial information, the GI/G/1 queue becomes a data-driven model that adjusts to available training data, for which this paper presents tight performance guarantees. 



The key idea of using MAD instead of variance as dispersion measure, is likely applicable to many other queueing system. 
Examples are queues with dependency and correlation structures in the series $ \{{\UU}_n \} $ and $ \{{\VV}_n \} $, the multi-server GI/G/$c$ queue and networks of queues. 
Indeed, most of the key performance measures for such systems are expectations of functions that are convex in the random variables (see e.g.,~\cite{shaked1988stochastic}), and therefore the mean-MAD approach can be used. 
The MAD perspective is of interest beyond queueing theory, because the search for extremal distributions of convex functions is relevant in many other settings. Moreover, whenever a performance measure can be viewed as a convex function of i.i.d.~random variables with mean-MAD ambiguity (e.g.,~nested max-operators in production systems; see \cite{glasserman1997bounds,bradley2002managing}), our approach will identify the extremal distribution and tight bounds.


The MAD approach stays close to the common practice in the stochastic field, namely to use probability distributions to model uncertainty. The nucleus of the MAD approach consists of the explicitly solvable dual LP described in Section~2. A simple reasoning then showed that this solution is independent of the precise objective function (in this paper describing waiting time moments of the GI/G/1 queue). Hence, the MAD approach is a generic, computationally tractable way to analyze stochastic processes, such as random walks and queues. 

Let us conclude with a broader robust optimization perspective. It is well-known that the use of probability distributions in stochastic systems often leads to computationally intractability (e.g., calculation of high dimensional convolutions).
Therefore, \cite{Bandi2012,bandi2015robust,whitt2017using} suggest to use uncertainty sets instead of probability distributions. The MAD approach described in this paper can serve in many situations as an alternative (not per se better), bringing new opportunities. The uncertainty set approach yields a worst-case scenario. Our approach yields both worst-case and best-case distributions, i.e., both upper and lower bounds. In stochastic systems one often studies convex functions in the stochastic variables. In the uncertainty set approach it is in general hard (in fact, NP-hard) to find worst-case scenarios for such convex functions. Our approach can easily find worst-case distributions as shown in this paper.

\ACKNOWLEDGMENT{The authors would like to thank Daniel Kuhn and Krzysztof Postek for pointing out the primal-dual reasoning that gives the intuitive proof of Theorem~\ref{thm1h}, and Marko Boon for helping with the experiments in Section~\ref{ec:contour}. 
}

\bibliographystyle{informs2014}
\bibliography{bibbook}

\ECSwitch

\pagestyle{fancy}
\renewcommand\headrulewidth{0pt}
\lhead{}\chead{}\rhead{}
\cfoot{\vspace*{1.5\baselineskip}\thepage}


\ECHead{E-Companion to ``MAD dispersion measure makes extremal queue analysis simple''}
\addtocontents{toc}{\setcounter{tocdepth}{-1}}

\section{Properties of MAD}\label{ap:properties_MAD}
We recall some well known properties of the MAD, see e.g.~\cite{BenTal1985}. Denote by $\sigma^2$ the variance of the random variable $X$, whose distribution is known to belong to the set $\mathcal{P}_{(\mu,d)}$. Then
$$
\frac{d^2}{4\beta (1-\beta)} \leq \sigma^2 \leq \frac{d(b-a)}{2}.
$$
In particular, since
$$
d^2 \leq 4\beta (1-\beta) \sigma^2 \leq \sigma^2,
$$
it holds that $d \leq \sigma$. For a proof, we refer the reader to \cite{BenTal1985}. For some distributions, an explicit formula for $d$ is available:
\begin{itemize}
\item Uniform distribution on $[a,b]$:
$$
d = \frac{1}{4}(b-a)
$$
\item Normal distribution $N(\mu,\sigma^2)$:
$$
d = \sqrt{\frac{2}{\pi}} \sigma
$$
\item Gamma distribution with parameters $\lambda$ and $k$ (for which $\mu = k/ \lambda$):
$$
d = \frac{2k^k}{\Gamma(k)\exp(k)} \frac{1}{\lambda}.
$$
\end{itemize}
The MAD is known to satisfy the bound
\begin{equation}
\label{eq:d_bound}
0 \leq d\leq \frac{2(b-\mu)(\mu-a)}{b-a}.
\end{equation}
Let $\beta=\mathbb{P}(X \geq \mu)$. For example, in the case of continuous symmetric distribution of $X$ we know that $\beta = 0.5$. This quantity is known to satisfy the bounds:
\begin{equation}
\label{eq:beta_bound}
\frac{d}{2(b-\mu)} \leq \beta \leq  1-\frac{d}{2(\mu - a)}.
\end{equation}

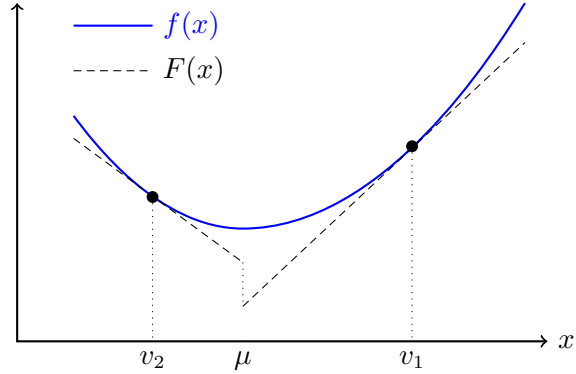
\begin{figure}[h!]
\begin{center}
\begin{tikzpicture}[scale=1.5]
\draw [<->,thick] (0,3) node (yaxis) [above] {}
        |- (4.7,0) node (xaxis) [right] {$x$};
\draw[blue,thick] (0.5,2) parabola bend (2,1) (4.5,3)
        node[below right] {};
        \draw[blue,thick] (0.5,2.8) -- (1.2,2.8) node[right] {$f(x)$};
        \draw[densely dashed] (0.5,2.4) -- (1.2,2.4) node[right] {$F(x)$};
        \draw[densely dashed] (2,0) -- (2,0) node[below] {$\mu$};
        \coordinate (a) at (0.5,1.8);
        \coordinate (f) at (2, 0.31);
        \coordinate (b) at (2,0.7);
        \coordinate (c) at (4.5,2.65);
         \draw[densely dashed] (a) -- (b);
         \draw[densely dashed] (f) -- (c);
        \coordinate (d) at (3.5,1.73);
        \fill[black] (d) circle (1.5pt);
        \coordinate (e) at (1.2,1.28);
        \fill[black] (e) circle (1.5pt);
         \draw[dotted] (b) -- (f);
          \draw[dotted] (e) -- (1.2,0) node[below] {$v_2$};
           \draw[dotted] (d) -- (3.5,0) node[below] {$v_1$};
\end{tikzpicture}
\caption{Some convex function $f(x)$ and its non-continuous piecewise linear minorant $F(x)$.}\label{fig:minor2}
\end{center}
\end{figure}

{
\section{Primal-dual proof of Theorem 3}\label{newproof3} 
In a similar manner as for the upper bound, we will show that the best-case distribution is a two-point distribution. We again consider the convex univariate measurable function $f(x)$ that has finite values on $[a,b]$. Under $\mathcal{P}_{(\mu,d,\beta)}$ ambiguity of the random variable $X$ we now need to solve
\begin{equation}\label{eq:primallb}
\begin{aligned}
&\min_{\Dist(x)\geq0} &  &\int_x f(x){\rm d} \Dist(x)\\
&\text{s.t.} &      & \int_x \mathbbm{1}_{\{x\geq\mu\}}{\rm d}\Dist(x)=\beta, \ \int_x |x-\mu|{\rm d}\Dist(x)=d, \ \int_x x{\rm d}\Dist(x)=\mu,\ \int_x {\rm d}\Dist(x)=1,   
\end{aligned}
\end{equation}
 which is a semi-infinite linear program with four equality constraints.


Consider the dual of \eqref{eq:primallb},
\begin{equation}\label{eq:duallb}
\begin{aligned}
&\max_{\lambda_0,\lambda_1,\lambda_2, \lambda_3} &  &\lambda_0 \beta + \lambda_1 d+\lambda_2 \mu+\lambda_3\\
&\text{s.t.} &      & f(x)-\lambda_0 \mathbbm{1}_{\{x\geq\mu\}} -\lambda_1|x-\mu|-\lambda_2x-\lambda_3\geq 0, \ \forall x\in[a,b].
\end{aligned}
\end{equation}
Define $F(x)=\lambda_0 \mathbbm{1}_{\{x\geq\mu\}} + \lambda_1|x-\mu|+\lambda_2x+\lambda_3$. Then the inequality in \eqref{eq:duallb} can be written as $F(x)\leq f(x)$, $\forall x$, i.e.~$F(x)$ minorizes $f(x)$. Note that in our new situation $F(x)$ has both a kink and a discontinuity at $x=\mu$, as depicted in Figure~\ref{fig:minor2}. The dual problem boils down to finding the tightest minorant that maximizes the dual problem's objective value. The minorant $F(x)$ touches the epigraph of $f(x)$ in at most two points on opposite sides of $\mu$ (i.e., $x_1\leq\mu\leq x_2$). This is a consequence of the supporting hyperplane theorem and the jump discontinuity at $x=\mu$. The dual problem now becomes
\begin{equation}\label{eq:duallb2}
\begin{aligned}
&\max_{\lambda_0,\lambda_1,\lambda_2, \lambda_3} &  &\lambda_0 \beta + \lambda_1 d+\lambda_2 \mu+\lambda_3\\
&\text{s.t.} &      & \lambda_0 +\lambda_1(x_1-\mu)+\lambda_2x_1+\lambda_3= f(x_1),  \\
& &      & \quad  -\lambda_1(x_2-\mu)+\lambda_2x_2+\lambda_3= f(x_2).
\end{aligned}
\end{equation}

Now using Lagrange duality, we can show that the optimal solution satisfies
$$
x_1=\mu+\frac{d}{2\beta}, \quad x_2=\mu-\frac{d}{2(1-\beta)},
$$
which corresponds to the values of $v_1$ and $v_2$ stated in Theorem~\ref{thm1hj}. Substituting this solution and solving for $\lambda_0, \lambda_1, \lambda_2$, and $\lambda_3$ gives
$$
\lambda_0=f(v_1)-f(v_2)+\frac{\lambda_1d}{(1-\beta)}-\frac{(\lambda_1+\lambda_2)d}{2\beta(1-\beta)}, \quad \lambda_3=f(v_2)+\frac{(\lambda_2 - \lambda_1)d}{2(1-\beta)}-\lambda_2\mu,
$$
and hence the objective value of the dual becomes $\beta f(v_1)+(1-\beta)f(v_2)$. Note that we have two free variables that can be chosen in a way that makes the solution dual feasible. The optimal probabilities of \eqref{eq:primallb} are obtained by solving the linear system resulting from \eqref{eq:primallb}, which produces the solution stated in Theorem~\ref{thm1hj}. Finally, one can verify that the primal and dual objective values are the same and that these results can be extended to the multivariate case in a manner analogous to that of Theorem~\ref{thm1h}.
}

\section{Representations for the tight bounds}\label{ec:contour}

We will now present some efficient ways of calculating the tight bounds identified in this paper. But first we show a way to verify the contour integral representation. 


\subsection{Numerical experiments with contour integrals}
Numerical aspects of integrals of the type \eqref{polgen2}
have been discussed in e.g.,~\cite{Abate1993,Janssen2015,chen2019extremal2}. 
For distributions with support on a finite set of points, 
potential numerical problems can arise, because $|{\rm Re}(\phi_X(u))|$ does not converge to zero as $|u|\to\infty$; see \cite{abate1992fourier,chen2019extremal2}. 
For the three-point distributions required in this paper we have performed extensive numerical experiments with \eqref{pollaczekup}. These experiments confirmed that the integrals can be calculated up to high accuracy with standard integration routines in Mathematica (our code is available upon request). 

For many parameter values $a,b,\mu,d$ such that \eqref{eq:d_bound} holds,
we have calculated $\expect{M}$ for generic increment $X_{(3)}$ using \eqref{pollaczekup}, and compared this with results from extensive stochastic simulations. We also compared the results with a third numerical procedure, known to be extremely stable and accurate. Let us explain the third procedure, which might be of independent interest. 
 
Choose the boundaries of the support as multiples of $\beta=|\mu|$ by writing that $a=-s\beta$ and $b=m\beta$ with $s,m$ positive integers. Denote by $M_\beta=M/\beta$ the normalized steady-state waiting time. We then get
$$
M_{\beta} \stackrel{d}{=}(M_\beta +X_\beta)^+,
$$
with $X_\beta=X/\beta$ a discrete random variable with support $\{-s,-1,m\}$ and MAD $$
d_\beta:=\expect{|X_\beta-  \expect{X_\beta}|}=
\frac{1}{\beta}\expect{|X-  \expect{X}|}=d.$$
Define $X_\beta=A_\beta-s$, so that
$$
M_{\beta} \stackrel{d}{=}(M_\beta +A_\beta-s)^+
$$
for a discrete random variable $A_\beta$ with support $\{0,s-1,s+m\}$ and probability generating function
$$
\expect{z^{A_\beta}}=p_a+p_\mu z^{s-1}+p_bz^{m+s},
$$
with
$$
p_a=\frac{d_\beta}{2(s-1)},\quad p_\mu=1-\frac{d_\beta}{2(s-1)}-\frac{d_\beta}{2(m+1)},\quad p_b=\frac{d_\beta}{2(m+1)}.
$$
Notice that $\expect{A_\beta}=s-1$. The resulting discrete queueing system is sometimes referred to as a bulk service queue. 
Let $r_0$ be the unique zero of $z^{s}-\expect{z^{A_\beta}}$ with real $z>1$. 
For any $\eps>0$ with $1+\eps<r_0$,
\beq \label{e111}
\expect{w^{M_\beta}}=\exp\Bigl(\frac{1}{2\pi i}\,\il_{|z|=1+\eps}\,{\rm ln}\Bigl(\frac{w-z}{1-z}\Bigr)\,\frac{(z^s-\expect{z^{A_\beta}})'}{z^s-\expect{z^{A_\beta}}}\,dz\Bigr)
\eq
holds when $|w|<1+\eps$. Alternatively,
\beq \label{e113}
\expect{w^{M_\beta}}=\frac{(s-\expect{A_\beta})(w-1)}{w^s-A(w)}\,\prod_{k=1}^{s-1}\,\frac{w-z_k}{1-z_k}
\eq
that holds for all $w$, $|w|<r_0$, in which $z_1,\ldots,z_{s-1}$ are the $s-1$ zeros of $z^s-\expect{z^{A_\beta}}$ in $|z|<1$.
Upon differentiation,  \eqref{e111} and  \eqref{e113} provide expressions for all cumulants of $M_\beta$ that are known to allow for accurate numerical evaluation, see \cite{Janssen2015}. 
We have then performed for a wide range of parameters, the following experiment:
\begin{henumerate}
\item Fix $\beta$, and then choose integers $s$ and $m$. In this way we create a standard bulk service queue with discrete-valued generic increment $A_\beta$. 
\item For ranging $d_\beta$, calculate $\expect{M_\beta}$ using root-finding procedures and  \eqref{e113}  or  using the contour integral  \eqref{e111}. 
\item Calculate
\begin{equation*}
\expect{M}=\frac{-1}{2\pi i}\int_{\mathcal{C}}\frac{\log(1-(p_a  {\rm e}^{-ua}+ p_b {\rm e}^{-ub}+p_c{\rm e}^{-uc}))}{u^2}{\rm d}u.
\end{equation*}
\item Check whether $\expect{M}=\beta \expect{M_\beta}$.
\end{henumerate}

\subsection{Numerical procedures for the GI/G/1 queue}
Calculations for $\expect{W_n}$ and $c_n(W)$ in the GI/G/1 queue can be performed using similar expressions as for the random walk. Let the random variable $\VV_{(3)}$ follow a three-point distribution on values $\{s_1,s_2,s_3\}$ with probabilities 
\begin{equation}\label{probs}
p_1=\frac{d_{\VV}}{2(\mu_{\VV}-a_{\VV})},\quad p_2=1-\frac{d_{\VV}}{2(\mu_{\VV}-a_{\VV})}-\frac{d_{\VV}}{2(b_{\VV}-\mu_{\VV})},\quad p_3=\frac{d_{\VV}}{2(b_{\VV}-\mu_{\VV})},
\end{equation}
with $0\leq a_{\VV}<\mu_{\VV}<b_{\VV}$, so that $\VV_{(3)}$ has mean $\mu_{\VV}$ and MAD $d_{\VV}$. Similarly, let $\UU_{(3)}$
have a three-point distribution on values $\{t_1,t_2,t_3\}$ with probabilities 
\begin{equation}\label{probs}
r_1=\frac{d_{\UU}}{2(\mu_{\UU}-a_{\UU})},\quad r_2=1-\frac{d_{\UU}}{2(\mu_{\UU}-a_{\UU})}-\frac{d_{\UU}}{2(b_{\UU}-\mu_{\UU})},\quad r_3=\frac{d_{\UU}}{2(b_{\UU}-\mu_{\UU})}
\end{equation}
and $0\leq a_{\UU}<\mu_{\UU}<b_{\UU}$, so that $\UU_{(3)}$ has mean $\mu_{\UU}$ and MAD $d_{\UU}$. 

We then have the representation, see also \cite{chen2019extremal},
\begin{equation}\label{permm13}
\expect{W_n}=
\sum_{k=1}^n \frac1k \sum_{\sum_i k_i=k,\sum_j l_j=k}\max\{0,\sum_{i=1}^3k_i s_i-\sum_{j=1}^3l_i t_i\}\cdot P(k_1,k_2,k_3)\cdot R(l_1,l_2,l_3)
\end{equation}
with 
$$
P(k_1,k_2,k_3)=\frac{k!}{k_1!k_2!k_3!}p_1^{k_1}p_2^{k_2}p_3^{k_3}, \quad R(l_1,l_2,l_3)=\frac{k!}{l_1!l_2!l_3!}r_1^{l_1}r_2^{l_2}r_3^{l_3},
$$
which requires summing $O(n^5)$ terms.

Let $\phi_{\VV_{(3)}}(s)$ and $\phi_{\UU_{(3)}}(s)$ denote the moment generating functions of $\VV_{(3)}$ and $\UU_{(3)}$. 
The tight upper bounds on $c_m(W)$ are given by
\begin{equation}\label{polgg1}
c_m(W)\leq \frac{(-1)^m}{2\pi i}\int_{\mathcal{C}}\frac{\log(1-\phi_{\VV_{(3)}}(-u)\phi_{\UU_{(3)}}(u))}{u^{m+1}}{\rm d}u,
\end{equation}
where $\mathcal{C}$ is a contour to the left of, and parallel to, the imaginary axis, and to the right of any singularities of $\log(1-\phi_{\VV_{(3)}}(-u)\phi_{\UU_{(3)}}(u))$ in the left half plane. Again comparing with extensive simulation, we have found the expression \eqref{polgg1} accurate and hence suitable for calculating the tight bounds.

\section{Distribution-free upper bounds for the GI/G/1 queue}\label{secdfb}
Consider the steady-state queue length $W$ in the GI/G/1 queue, which satisfies $W\das (W+\VV-\UU)^+$. Denote by $\sigma^2_\UU$ and $\sigma^2_\VV$ the variances of $\UU$ and $\VV$, respectively. Let $\rho=\expect{\VV}/\expect{\UU}<1$. The following bounds on $\expect{W}$ only require information about the first two moments of $\UU$ and $\VV$:
\begin{itemize}
\item Kingman's upper bound:
\begin{equation}\label{kingman}
\expect{W}\leq\frac{\sigma_\VV^2+\sigma_\UU^2}{2(\expect{\UU}-\expect{\VV})}.
\end{equation}
\item Daley's upper bound:
\begin{equation}\label{daley}
\expect{W}\leq\frac{\sigma_\VV^2+\rho(2-\rho)\sigma_\UU^2}{2(\expect{\UU}-\expect{\VV})}.
\end{equation}
\item 
Upper bound of \cite{chen2019extremal} based on the two-point conjecture:
\begin{equation}\label{chenwhittbound}
\expect{W}\leq\frac{\sigma_\VV^2+\kappa(\rho)\sigma_\UU^2}{2(\expect{\UU}-\expect{\VV})},
\end{equation}
with $\kappa(\rho)=2\rho(1-\rho)/(1-\delta)$ and $\delta\in(0,1)$ the solution of $\delta=\exp(-(1-\delta)/\rho)$. 
\end{itemize}

{
\section{Further numerical results for the bounds}\label{sec:ECn}
We now complement Table \ref{tab11} with some more numerical values for the bounds on $\expect{W}$. Table \ref{tab11x} gives the unscaled values of $\expect{W}$ for the same parameter values as in Table \ref{tab11}. 

\begin{table}[h]

\caption{{Bounds for $\expect{W}$  for $(\mu_{\UU},d_{\UU},a_{\UU},b_{\UU})=(1,1,0,10)$ and 
 $(\mu_{\VV},d_{\VV},a_{\VV},b_{\VV})=(\rho,0.1,0,10)$.}}\label{tab11x}
 \vspace{.3cm}
\begin{center}

\begin{tabular}{r R{1.8cm} R{1.8cm} R{1.8cm} R{1.8cm}}\hline 
   $\rho$    & Tight (Thm.~\ref{tightgg1}) & C \& W \eqref{chenwhittbound} & Daley \eqref{daley} & Kingman \eqref{kingman}  \\ \hline
 0.1 & 0.45179 & 0.77780 & 0.80556 & 3.05556 \\
 0.2 & 0.63077 & 1.31953 & 1.43750 & 3.43750 \\
 0.5 & 2.03141 & 3.63750 & 4.25000 & 5.50000 \\
 0.7 & 5.81373 & 7.39989 & 8.41667 & 9.16667 \\
 0.8 & 10.47728 & 12.02090 & 13.25000 & 13.75000 \\
 0.9 & 24.28220 & 25.80400 & 27.25000 & 27.50000 \\
 0.95 & 51.79564 & 53.31910 & 54.87500 & 55.00000 \\
 0.99 &  271.80153 & 273.33100 & 274.97500 & 275.00000 \\  \hline
\end{tabular}
\end{center}
\end{table}

The variance bounds are often reported in terms of the squared coefficient of variation (variance divided by the square of the mean), see  \cite{chen2019extremal}. For the extremal distributions with $(\mu_{\VV},d_{\VV},a_{\VV},b_{\VV})=(\rho,d_V,0,b_V)$
and 
$(\mu_{\UU},d_{\UU},a_{\UU},b_{\UU})=(1,d_U,0,b_U)$
this gives
\begin{equation*}
c_V^2=\frac{\sigma^2_V}{\mu_V^2}=\frac{d_Vb_V}{2\rho^2}, \quad c_U^2=\frac{\sigma^2_U}{\mu_U^2}=\frac{d_Ub_U}{2}.
\end{equation*}

Fixing the squared coefficient of variations $c_V^2$ and $c_U^2$ is equivalent with choosing the MADs as 
\begin{equation}\label{eqcov}
d_V=\frac{2\rho^2c_V^2}{b_V}, \quad d_U=\frac{2 c_U^2}{b_U}.
\end{equation}

We next present in Tables \ref{tab2x}-\ref{tab2x44} some further numerical results, for $c_U^2=c_V^2=0.5$ and $c_U^2=c_V^2=4$. 

\begin{table}[h!]
\caption{Bounds for 
$\expect{W}$  for  $(\mu_{\VV},d_{\VV},a_{\VV},b_{\VV})=(\rho,d_{\VV},0,10)$ and $(\mu_{\UU},d_{\UU},a_{\UU},b_{\UU})=(1,d_{\UU},0,10)$ 
 with $d_{\VV}, d_{\UU}$ as in \eqref{eqcov} and $c_U^2=c_V^2=0.5$.}\label{tab2x}
 \vspace{.3cm}
\begin{center}
\begin{tabular}{r R{1.5cm} R{1.5cm} R{1.5cm} R{1.5cm}}\hline
   $\rho$    & Tight (Thm.~\ref{tightgg1}) & C \& W \eqref{chenwhittbound} & Daley \eqref{daley} & Kingman \eqref{kingman}  \\ \hline
 0.1 & 0.00785 & 0.05278 & 0.05555 & 0.28055 \\
 0.2 & 0.02230 & 0.11320 & 0.12500 & 0.32500 \\
 0.5 & 0.14921 & 0.43875 & 0.50000 & 0.62500 \\
 0.7 & 0.48818 & 1.06499 & 1.16667 & 1.24167 \\
 0.8 & 0.99509 & 1.87709 & 2.00000 & 2.05000 \\
 0.9 & 2.85149 & 4.35540 & 4.50000 & 4.52500 \\
 0.95 & 7.29378& 9.34441 & 9.50000 & 9.51250 \\
 0.99 & 46.78335 & 49.33560 & 49.50000 & 49.50250 \\ \hline
\end{tabular}
\end{center}
\end{table}

\begin{table}[h!]
\caption{Bounds for 
$\expect{W}$  for  $(\mu_{\VV},d_{\VV},a_{\VV},b_{\VV})=(\rho,d_{\VV},0,10)$ and $(\mu_{\UU},d_{\UU},a_{\UU},b_{\UU})=(1,d_{\UU},0,10)$ 
 with $d_{\VV}, d_{\UU}$ as in \eqref{eqcov} and $c_U^2=c_V^2=4$.}\label{tab2x4}
 \vspace{.3cm}
\begin{center}
\begin{tabular}{r R{1.8cm} R{1.8cm} R{1.8cm} R{1.8cm}}\hline
   $\rho$    & Tight (Thm.~\ref{tightgg1}) & C \& W \eqref{chenwhittbound} & Daley \eqref{daley} & Kingman \eqref{kingman}  \\ \hline
 0.1 & 0.09358 & 0.42224 & 0.44444 & 2.24444 \\
 0.2 & 0.26429 & 0.90562 & 1.00000 & 2.60000 \\
 0.5 & 2.05142 & 3.51000 & 4.00000 & 5.00000 \\
 0.7 & 6.76335 & 8.51991 & 9.33333 & 9.93333 \\
 0.8 & 13.18168 & 15.01670 & 16.00000 & 16.40000 \\
 0.9 & 32.95685 & 34.84320 & 36.00000 & 36.20000 \\
 0.95 & 72.84232 & 74.75520 & 76.00000 & 76.10000 \\
 0.99 & 392.74278 & 394.68400 & 396.00000 & 396.02000 \\ \hline
\end{tabular}
\end{center}
\end{table}

\begin{table}[h!]
\caption{{Bounds for 
$\expect{W}$  for  $(\mu_{\VV},d_{\VV},a_{\VV},b_{\VV})=(\rho,d_{\VV},0,10)$ and $(\mu_{\UU},d_{\UU},a_{\UU},b_{\UU})=(1,d_{\UU},0,10)$ 
 with $d_{\VV}, d_{\UU}$ as in \eqref{eqcov}, $c_U^2=4$ and $c_V^2=0.5$.}}\label{tab2x0}
 \vspace{.3cm}
\begin{center}
\begin{tabular}{r R{1.8cm} R{1.8cm} R{1.8cm} R{1.8cm}}\hline
   $\rho$    & Tight (Thm.~\ref{tightgg1}) & C \& W \eqref{chenwhittbound} & Daley \eqref{daley} & Kingman \eqref{kingman}  \\ \hline
0.1  & 0.07003   & 0.40280   & 0.42500   & 2.22500   \\
0.2  & 0.15280   & 0.81812   & 0.91250   & 2.51250   \\
0.5  & 0.91273   & 2.63500   & 3.12500   & 4.12500   \\
0.7  & 3.73777   & 5.66158   & 6.47500   & 7.07500   \\
0.8  & 7.53710   & 9.41674   & 10.40000  & 10.80000  \\
0.9  & 18.82048  & 20.66820  & 21.82500  & 22.02500  \\
0.95 & 41.31986  & 43.16770  & 44.41250  & 44.51250  \\
0.99 & 221.30939 & 223.16700 & 224.48200 & 224.50200 \\ \hline
\end{tabular}
\end{center}
\end{table}

\begin{table}[h!]
\caption{{Bounds for 
$\expect{W}$  for  $(\mu_{\VV},d_{\VV},a_{\VV},b_{\VV})=(\rho,d_{\VV},0,10)$ and $(\mu_{\UU},d_{\UU},a_{\UU},b_{\UU})=(1,d_{\UU},0,10)$ 
 with $d_{\VV}, d_{\UU}$ as in \eqref{eqcov}, $c_U^2=0.5$ and $c_V^2=4$.}}\label{tab0x2}
 \vspace{.3cm}
\begin{center}
\begin{tabular}{r R{1.8cm} R{1.8cm} R{1.8cm} R{1.8cm}}\hline
   $\rho$    & Tight (Thm.~\ref{tightgg1}) & C \& W \eqref{chenwhittbound} & Daley \eqref{daley} & Kingman \eqref{kingman}  \\ \hline
0.1  & 0.02599   & 0.07222   & 0.07500   & 0.30000   \\
0.2  & 0.10463   & 0.20070   & 0.21250   & 0.41250   \\
0.5  & 1.00498   & 1.31375   & 1.37500   & 1.50000   \\
0.7  & 3.39670   & 3.92332   & 4.02500   & 4.10000   \\
0.8  & 6.81534   & 7.47709   & 7.60000   & 7.65000   \\
0.9  & 17.72431  & 18.53040  & 18.67500  & 18.70000  \\
0.95 & 40.05188  & 40.93190  & 41.08750  & 41.10000  \\
0.99 & 219.91292 & 220.85300 & 221.01700 & 221.02000 \\ \hline
\end{tabular}
\end{center}
\end{table}

\begin{table}[h!]
\caption{Bounds for 
$(1-\rho)\expect{W}/\rho$  for  $(\mu_{\VV},d_{\VV},a_{\VV},b_{\VV})=(\rho,d_{\VV},0,10)$ and $(\mu_{\UU},d_{\UU},a_{\UU},b_{\UU})=(1,d_{\UU},0,10)$ 
 with $d_{\VV}, d_{\UU}$ as in \eqref{eqcov} and $c_U^2=c_V^2=0.5$.}\label{tab2xx}
 \vspace{.3cm}
\begin{center}
\begin{tabular}{r R{1.5cm} R{1.5cm} R{1.5cm} R{1.5cm}}\hline
   $\rho$    & Tight (Thm.~\ref{tightgg1}) & C \& W \eqref{chenwhittbound} & Daley \eqref{daley} & Kingman \eqref{kingman}  \\ \hline
 0.1 & 0.07070 & 0.47502 & 0.50000 & 2.52500 \\
 0.2 & 0.08922 & 0.45281 & 0.50000 & 1.30000 \\
 0.5 & 0.14921 & 0.43875 & 0.50000 & 0.62500 \\
 0.7 & 0.20922 & 0.45642 & 0.50000 & 0.53214 \\
 0.8 & 0.24877 & 0.46927 & 0.50000 & 0.51250 \\
 0.9 & 0.31683 & 0.48393 & 0.50000 & 0.50277 \\
 0.95 & 0.38388 & 0.49181 & 0.50000 & 0.50065 \\
 0.99 & 0.47255 & 0.49833 & 0.50000 & 0.50002 \\ \hline
\end{tabular}
\end{center}
\end{table}

\begin{table}[h!]
\caption{Bounds for 
$(1-\rho)\expect{W}/\rho$  for  $(\mu_{\VV},d_{\VV},a_{\VV},b_{\VV})=(\rho,d_{\VV},0,10)$ and $(\mu_{\UU},d_{\UU},a_{\UU},b_{\UU})=(1,d_{\UU},0,10)$ 
 with $d_{\VV}, d_{\UU}$ as in \eqref{eqcov} and $c_U^2=c_V^2=4$.}\label{tab2x44}
 \vspace{.3cm}
\begin{center}
\begin{tabular}{r R{1.5cm} R{1.5cm} R{1.5cm} R{1.5cm}}\hline
   $\rho$    & Tight (Thm.~\ref{tightgg1}) & C \& W \eqref{chenwhittbound} & Daley \eqref{daley} & Kingman \eqref{kingman}  \\ \hline
 0.1 & 0.84228 & 3.80016 & 4.00000 & 20.20000 \\
 0.2 & 1.05719 & 3.62248 & 4.00000 & 10.40000 \\
 0.5 & 2.05142 & 3.51000 & 4.00000 & 5.00000 \\
 0.7 & 2.89858 & 3.65139 & 4.00000 & 4.25714 \\
 0.8 & 3.29542 & 3.75418 & 4.00000 & 4.10000 \\
 0.9 & 3.66187 & 3.87146 & 4.00000 & 4.02222 \\
 0.95 & 3.83381 & 3.93449 & 4.00000 & 4.00526 \\
 0.99 & 3.96710 & 3.98671 & 4.00000 & 4.00020 \\ \hline
\end{tabular}
\end{center}
\end{table}

\begin{table}[h!]
\caption{{Bounds for 
$(1-\rho)\expect{W}/\rho$   for  $(\mu_{\VV},d_{\VV},a_{\VV},b_{\VV})=(\rho,d_{\VV},0,10)$ and $(\mu_{\UU},d_{\UU},a_{\UU},b_{\UU})=(1,d_{\UU},0,10)$ 
 with $d_{\VV}, d_{\UU}$ as in \eqref{eqcov}, $c_U^2=4$ and $c_V^2=0.5$.}}\label{tab2x0x}
 \vspace{.3cm}
\begin{center}
\begin{tabular}{r R{1.8cm} R{1.8cm} R{1.8cm} R{1.8cm}}\hline
   $\rho$    & Tight (Thm.~\ref{tightgg1}) & C \& W \eqref{chenwhittbound} & Daley \eqref{daley} & Kingman \eqref{kingman}  \\ \hline
0.1  & 0.63030 & 3.62516 & 3.82500 & 20.02500 \\
0.2  & 0.61120 & 3.27248 & 3.65000 & 10.05000 \\
0.5  & 0.91273 & 2.63500 & 3.12500 & 4.12500  \\
0.7  & 1.60190 & 2.42639 & 2.77500 & 3.03214  \\
0.8  & 1.88427 & 2.35418 & 2.60000 & 2.70000  \\
0.9  & 2.09116 & 2.29646 & 2.42500 & 2.44722  \\
0.95 & 2.17473 & 2.27199 & 2.33750 & 2.34276  \\
0.99 & 2.23545 & 2.25421 & 2.26750 & 2.26770 \\ \hline
\end{tabular}
\end{center}
\end{table}

\begin{table}[h!]
\caption{{Bounds for 
$(1-\rho)\expect{W}/\rho$   for  $(\mu_{\VV},d_{\VV},a_{\VV},b_{\VV})=(\rho,d_{\VV},0,10)$ and $(\mu_{\UU},d_{\UU},a_{\UU},b_{\UU})=(1,d_{\UU},0,10)$ 
 with $d_{\VV}, d_{\UU}$ as in \eqref{eqcov}, $c_U^2=0.5$ and $c_V^2=4$.}}\label{tab0x2x}
 \vspace{.3cm}
\begin{center}
\begin{tabular}{r R{1.8cm} R{1.8cm} R{1.8cm} R{1.8cm}}\hline
   $\rho$    & Tight (Thm.~\ref{tightgg1}) & C \& W \eqref{chenwhittbound} & Daley \eqref{daley} & Kingman \eqref{kingman}  \\ \hline
0.1  & 0.23392 & 0.65002 & 0.67500 & 2.70000 \\
0.2  & 0.41852 & 0.80281 & 0.85000 & 1.65000 \\
0.5  & 1.00498 & 1.31375 & 1.37500 & 1.50000 \\
0.7  & 1.45573 & 1.68142 & 1.72500 & 1.75714 \\
0.8  & 1.70384 & 1.86927 & 1.90000 & 1.91250 \\
0.9  & 1.96937 & 2.05893 & 2.07500 & 2.07778 \\
0.95 & 2.10799 & 2.15431 & 2.16250 & 2.16316 \\
0.99 & 2.22134 & 2.23084 & 2.23250 & 2.23253 \\ \hline
\end{tabular}
\end{center}
\end{table}

}

\end{document}